\documentclass[a4paper,11pt]{amsart}
\usepackage{amsaddr}
\usepackage{amsfonts,amsopn}
\usepackage{rotating,graphicx}
\usepackage{nccrules}
\usepackage[left=2.5cm, right=2.5cm, width=16cm, height=25cm ]{geometry}

\usepackage{tikz,fp}

\newcommand{\bs}{\mathbf}
\newcommand{\pd}[2]{\dfrac{\partial #1}{\partial #2}}

\let\bm\boldsymbol

\begin{document}


\title[An Immersed method based on cut-cells for the simulation of 2D incompressible]{An Immersed method based on cut-cells for the simulation of 2D incompressible fluid flows past solid structures}
\author{Fran\c{c}ois Bouchon}
\address[A. One]{Laboratoire de Math\'ematiques Blaise Pascal, UMR 6620, \\
         Universit\'e Clermont Auvergne and CNRS, \\
         Campus des C\'ezeaux, 3 place Vasarely, TSA 60026 CS 60026, \\
         63178 Aubi\`ere cedex, France.}
\email[A. One]{francois.bouchon@uca.fr}

\author{Thierry Dubois}
\address[A. Two]{Laboratoire de Math\'ematiques Blaise Pascal, UMR 6620, \\
         Universit\'e Clermont Auvergne and CNRS, \\
         Campus des C\'ezeaux, 3 place Vasarely, TSA 60026 CS 60026, \\
         63178 Aubi\`ere cedex, France.}
\email[A. Two]{thierry.dubois@uca.fr}

\author{Nicolas James}
\address[A. Three]{Laboratoire de Math\'ematiques et Applications UMR 7348, \\
  Universit\'e de Poitiers, \\
  T\'el\'eport 2 - B.P. 30179, Boulevard Marie et Pierre Curie, \\
  86962 Futuroscope Chasseneuil Cedex, France.}
\email[A. Three]{nicolas.james@math.univ-poitiers.fr}

\begin{abstract}
We present a cut-cell method for the simulation of 2D incompressible flows past obstacles. It consists in using the MAC scheme on cartesian grids and imposing Dirchlet boundary conditions for the velocity field on the boundary of solid structures following the Shortley-Weller formulation. In order to ensure local conservation properties, viscous and convecting terms are discretized in a finite volume way. The scheme is second order implicit in time for the linear part, the linear systems are solved by the use of the capacitance matrix method for non-moving obstacles. Numerical results of flows around an impulsively started circular cylinder are presented which confirm the efficiency of the method, for Reynolds numbers 1000 and 3000. An example of flows around a moving rigid body at Reynolds number 800 is also shown, a solver using the PETSc-Library has been prefered in this context to solve the linear systems.
\end{abstract}


\maketitle

\section{Introduction}
\label{intro}
For some decades, many researchers and engineers have been considering the numerical solution of fluid flows, for different kind of fluids and different geometries. With the increasing performance of super-computers, it has been possible to tackle more and more challenging problems, for higher Reynolds numbers and complex geometries. Several different discretization techniques can be used to consider these problems: Finite Element Methods, Finite Volumes Methods, Spectral Methods, and Finite Difference Methods.
The MAC scheme on cartesian grids \cite{Harlow_Welch} can be viewed both as a Finite Volume Method or a Finite Difference Method on staggered grids, and is adapted to 2D or 3D flows in simple geometries, for example for lid-driven cavity or backward facing step.
To take into account some obstacles in the flows (or complex geometries), immersed boundary techniques have been developped by Peskin in the 80's (\cite{Peskin,Peskin2}), consisting in using Dirac functions to model the interacting force between the fluid and the solid structure. These methods have inspired many authors in the following years, Mohd-Yusof has combined them with the use of B-Splines (\cite{Mohd-Yusof}) in his momentum forcing methods to consider complex geometries. The main advantage of these techniques is that the forcing term does not change the spatial operators, making them quite easy to implement (see \cite{Mittal} for a review, and refences therein).
As an alternative, Bruno \textit{et al.} have developped penalization techniques to inforce suitable boundary conditions \cite{Bruneau_al}. Similar techniques have also been investigated by Maury \textit{et al.} \cite{Maury_al} and justified from a mathematical point of view in \cite{Maury}. These methods have been shown to be efficient in the context of several particles in a flow \cite{Lefebvre}, and when considering possible collisions between them \cite{Verdon_al}.

Arbitrary Lagrangian Eulerian (ALE) methods have been developped for flows in geometries which vary in time (see \cite{Richter,Richter2,Richter3} where authors use some of the ideas of \cite{Belytschko}). The aim is to formulate the equation in a fixed reference domain, by using a mapping from the reference domain $\Omega(0)$ to the domain $\Omega(t)$ occupied by the fluid at time $t$. The position of the moving bodies, which correspond to the boundary of $\Omega(t)$, being available, the velocity field of these bodies defined on $\partial\Omega(t)$ have to be extended to $\Omega(t)$. Once this is done (generally with harmonic extensions), the equations are written in the reference domain by using the chain-rule formula.

For problems involving non-rigid bodies, Roshchenko \textit{et al.} (see \cite{Minev}) have used splitting methods to solve first the evolution of the velocity field in the fluid, and then to consider the deformation of the body. These ideas of splitting the model can be viewed as similar to the projection techniques (see \cite{Cho}, and \cite{GMJ} for a review  and references therein).

The method presented in this work joints another family of methods, called cut-cell methods. The idea of these methods is to modify the discretization of the Navier-Stokes Equations in the cells cut by the immersed boundary (see  \cite{ye_mittal_udaykumar_shyy,Mittal_al,Tucker,Chung}). One can discretize the equations on smallest cells obtained by intersecting the grid-cell with the domain occupied by the fluid, or one can merge these smallest cells with a neighbouring one. These methods can be combined with the levelset methods to track the boundary of the fixed or moving body (see \cite{Sethian}).
These ideas are used in the present work: the body in the fluid is represented by a levelset function, and the location of the velocity components are modified in the cut-cells (see \cite{BDJ1}), the pressure remaining placed at the center of cartesian grid cells for both fluid-cells and cut-cells. For the Laplacian of the velocity, the classical five-point approximation must be replaced by a local 6-point formula, for which the truncation error is only first order. But as in \cite{ Matsu_Yama}, global second order convergence of the method is recovered.
This second-order convergence for the velocity and the pressure with our cut-cell scheme has been obtained for the flow past a circular cylinder at Re=40 in \cite{Nico} by comparing with the reference solution proposed in \cite{Lamballais}.

The paper is organized as follows: Section \ref{sec:setting} is devoted to the presentation of the problem. The Navier-Stokes Equations are considered in a 2D geometry, which is supposed to be fixed for the sake of clarity. We also introduce there the notation for the grids, and detail space discretization. In Section \ref{sec:comput}, we give some information about computational aspects.
In the case of fixed domain, a fast solver adapted from the capacitance matrix method (see \cite{Buzbee_Golub,Buzbee_Dorr}) is used to solve the linear systems for both components and for the pressure. We also show that the method can be adapted in the case of moving domain. In this context, the preprocessing step of the capacitance matrix method would need to be done at each time-iteration which would then increase the CPU time. Therefore, we have prefered for this case a parallel version of an algebraic multigrid method (HYPRE BoomerAMG) implemented using the PETSc Fortran library (see~\cite{petsc-web-page,petsc-user-ref}).

Section~\ref{sec:numer} is then devoted to numerical tests, where we compare our results with theoretical predictions and other numerical results available in the literature. A conclusion is given in section~\ref{sec:conclusion}.

\section{The Setting of the Problem}
\label{sec:setting}

\subsection{Flows past obstacles}
\begin{figure}
\includegraphics[width=0.6\textwidth]{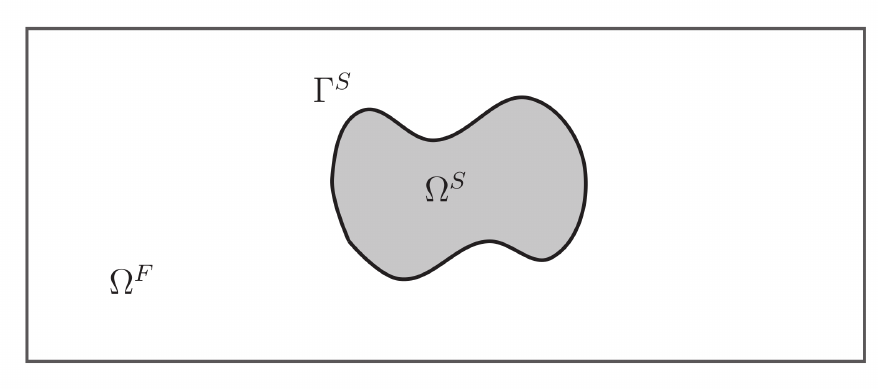}
\caption{The solid body $\Omega^S$ with boundary $\Gamma^S$ and the surrounding computational domain $\Omega^F$ in which the flow is to be simulated.}
\label{fig-domain}
\end{figure}

We consider a flow in a two-dimensional domain $\Omega=(0,L)\times(0,H)$ which contains a domain $\Omega^S$ occupied by the solid which is supposed to be fixed for the sake of simplicity. We denote then $\Omega^F=\Omega\setminus\Omega^S$ the domain occupied by the fluid (see Fig. \ref{fig-domain}).

The velocity field in the fluid satisfies the Navier-Stokes equations, with no-slip boundary conditions. We consider then the problem:

\begin{align}
  & \pd{\bs u}{t}\,-\,\frac1{Re}\,\Delta{\bs u}\,+\,\bs{\nabla}\cdot({\bs u}\otimes\bs{u})\,+\,{\bs\nabla}p\,=\,0, \label{nse_1} \\
  & {\bs\nabla}\cdot{\bs u}\,=\,0, \label{nse_2} \\
  & {\bs u}({\bs x},t=0)\,=\,{\bs u}_0,\label{nse_3}
\end{align}
where ${\bs u}({\bs x},t)=(u,v)$ is the velocity field at ${\bs x}=(x,y)\in\Omega^F$ at time $t>0$,
${\bs u}_0$ is the initial condition and $Re$ is the Reynolds number.
We impose homogeneous Dirichlet boundary conditions for the velocity field on $\partial{\Omega^F}$:
\begin{align}
& {\bs u}\,=\,0\mbox{ on }\partial\Omega^F\label{nse_bc}
\end{align}
We mention that non-homogenous Dirichlet boundary conditions can also be treated with the method presented here.
\subsection{Discretization}

For the time-discretization of \eqref{nse_1}--\eqref{nse_3}, we use a second-order backward difference (BDF2) 
projection scheme. In a first step, the velocity field is advanced in time with a semi-implicit scheme
decoupling the velocity and pressure unknowns. Then, the intermediate velocity is projected in order to obtain
a free-divergence velocity field.

Let $\delta t>0$ stand for the time step and $t^k=k\,\delta t$ discrete time values.
Let us consider that $(\bs u^j,\, P^j)$ are known for $j\le k$. The computation of $(\bs u^{k+1},\, P^{k+1})$ 
needs two steps:

\begin{equation}\label{nse_proj-1}
   \dfrac{3\widetilde{\bs u}^{k+1}-4{\bs u}^k+{\bs u}^{k-1}}{2\delta t}  - \frac 1{Re}\,\Delta\widetilde{\bs u}^{k+1} + \bs\nabla P^k =
    - 2\,\bs\nabla\cdot(\bs u^k\otimes\bs u^k)+\bs\nabla\cdot(\bs u^{k-1}\otimes\bs u^{k-1})  
\end{equation}
with homogeneous Dirichlet boundary condition for $\widetilde{\bs u}^{k+1}$.

Then the intermediate velocity field $\widetilde{\bs u}^{k+1}$ is projected in the free-divergence space to get $\bs{u}^{k+1}$:
\begin{equation}\label{nse_proj-2}
  \begin{split}
   &\dfrac{\bs{u}^{k+1}-\widetilde{\bs u}^{k+1}}{\delta t} + \frac 23\,\bs\nabla(P^{k+1}-P^k) = 0,  \\
   &\bs{\nabla}\cdot\bs{u}^{k+1} = 0, \quad (\bs{u}^{k+1}-\widetilde{\bs u}^{k+1})\cdot\bs{n} = 0 \quad \text{on }\Gamma.    
  \end{split}
\end{equation}

\begin{figure}
\includegraphics[width=0.6\textwidth]{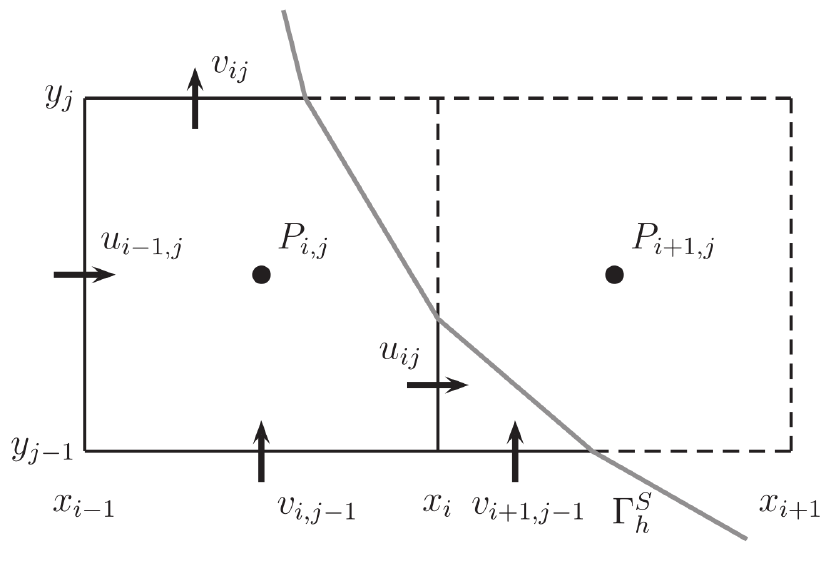}
\caption{Location of the unknowns near the solid body}
\label{fig-locationunknown}
\end{figure}

\begin{figure}
\includegraphics[width=0.6\textwidth]{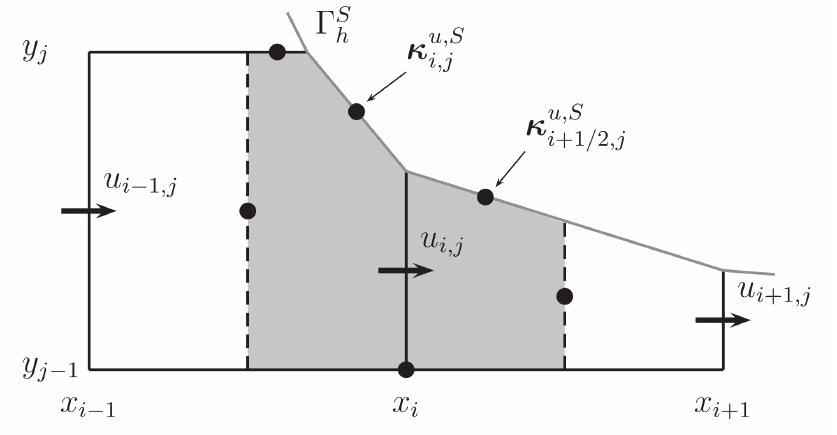}
\caption{Location of the computation of the fluxes near the solid body}
\label{fig-flux}
\end{figure}

For the spatial discretization, we modify the MAC scheme near the boundary by changing the location of the unknowns of the velocity components for the cells cut by the solid as depicted on Fig. \ref{fig-locationunknown} , the pressure unknowns remaining in their original place
(see \cite{BDJ1} for more details).
To discretize the Laplacian in \eqref{nse_proj-1}, we must replace the five-point formula by a six-point discretization. For the convective terms, the fluxes are computed at the midle of the vertical and horizontal edges (see Fig.  \ref{fig-flux}).
For the pressure, linear interpolation are used rather than changing the location of the pressure unknowns. The same kind of linear interpolation is used to get consistant evaluation of the pressure gradient in \eqref{nse_proj-2}.

Although the truncation error is only first order in space for the resulting numerical scheme, the second order accuracy is recovered which is due to a superconvergence phenomena analoguous to those proven in \cite{Matsu_Yama}.
This second order has been observed in \cite{Nico} by showing results in comparison with those of \cite{Lamballais}.

\section{Computational Aspects}
\label{sec:comput}
\subsection{A fast parallel direct solver to treat fixed solid strutures }
When considering fixed solid bodies, the use of a direct solver with a preprocessing procedure is efficient. Once the preprocessing computations have been done, the cost paid to solve the linear
systems is about twice the case of a numerical simulation in the same computational domain without obstacles. We summarize hereafter the fast direct solver derived from the capacitance matrix method and adapted to the case of non
uniform grids (see \cite{BDJ1} and \cite{BP2010} for the details) which has been implemented in our code.

After spatial discretization of the Navier-Stokes equations, one linear equation is obtained per node in the part of the computational domain
filled by the fluid and per unknown that is $u,\, v$ and $p$. We complete these sets of linear equations by adding similar ones for nodes of the cartesian grid lying inside
the solid obstacle but with zero as right-hand side. The unknowns corresponding to mesh points in $\Omega_h^S$ are fictitious ones. As in $\Omega^F$, the numerical scheme accounts for
the boundary conditions on $\Gamma_h^S$, the fluid unknowns are independant to the solid ones.
We therefore obtain linear algebraic systems defined on the whole cartesian grid with $n_x\times n_y$ mesh points whose sizes are $(n_x-1)\times n_y$ for $u$, $n_x\times (n_y-1)$ for $v$
and $(n_x-1)\times (n_y-1)$ for $p$. All three linear systems are similar in nature~: the resulting matrices have similar structures with five or six non-zero coefficients per row.

Let us consider one of these linear system. We denote by $N$ its size and by $A\in \mathcal{M}_N(\mathbb{R})$ its matrix. Then at each time iteration, we have to solve a linear system
\begin{equation} \label{Ax=z}
   AX\,=\, Z
\end{equation}
with the right-hand side $Z$ computed from the velocity and the pressure at previous time steps. As it is mentioned above, the matrix $A$ is non-symmetric. Let us consider now the matrix $G$ obtained
with the same discretization on the whole computational domain $\Omega$ totally filled by a fluid that is with no obstacles. The matrices $A$ and $G$ differ only on rows corresponding to
computational meshes for which the five-point stencil interacts with a cut-cell. Let us denote by $n_c$ this total number of rows, namely rows such that $A-G$ have
non-vanishing coefficients.
The efficiency of our direct solver is due to the fact that $n_c$ is small compared with $N$ and that the non-zero coefficients on each row of $A-G$ is bounded.
The linear system \eqref{Ax=z} can be rewritten as
\begin{equation} \label{Gx=z-Qy}
   GX\, =\, Z\,-\, QY
\end{equation}
where $Q$ is a matrix of dimensions $N\times n_c$ with one non-vanishing coefficient per column, equal to one, and $Y\in\mathbb{R}^{n_c}$ such that
\begin{equation*}
  Q\,Y\,=\, (A-G)\,X.
\end{equation*}
It can be easily shown, using $Q^tQ=I_{n_c}$, that $Y$ is solution of the following linear system
\begin{equation}  \label{Comp_y}
   \left( I_{n_c}\,+\, M\,G^{-1}Q\right)\, Y\,=\, M\, G^{-1}\,Z
\end{equation}
with $M=Q^t(A-G)$. The matrix $I_{n_c}\,+\, M\,G^{-1}Q$ is a non-singular matrix (see \cite{BP2010} for a proof) of size $n_c$.

Based on these relations, the algorithm implemented to solve \eqref{Ax=z} consists in a preprocessing step where the matrix $I_{n_c}\,+\, M\,G^{-1}Q$ is factorized
(we use a $LU$-factorization) followed by
\begin{enumerate}
   \item[i)]   Compute $Z$ and solve $GW=Z$~;
   \item[ii)]  Compute $MW$ and solve \eqref{Comp_y}~;
   \item[iii)] Compute $QY$ and solve $GX\,=\,Z\,-\,QY$.
\end{enumerate}
Recalling that $G$ is the matrix corresponding to the standard MAC scheme on the whole computational mesh, steps i) and iii) can be performed by using
any efficient solvers available on cartesian grids. In the present work, we use Discrete Fourier transforms in the vertical direction (where the mesh
is uniform) combined with $LU$-factorizations of the resulting tridiagonal systems.

The parallel version of this direct solver is based on explicit communications performed by calling functions of the MPI library. The main feature
of MPI is that a parallel application consists in running $p$ independant processes which may be executed on different computers, processors or
cores. These processes can exchange datas by sending/receiving messages \textit{via} a network connecting all the involved computing units.

\begin{figure}
  \begin{center}
\begin{tikzpicture}[scale=0.7]
\def\ymin{7}
\def\Height{4}
\FPadd\ymax{\ymin}{\Height}

\def\halflength{1}

\draw[thick,dashed] (2,\ymin)  -- (4,\ymin) ;
\draw[thick]        (4,\ymin)  -- (12,\ymin) ;
\draw[thick,dashed] (12,\ymin) -- (14,\ymin) ;

\draw[thick,dashed] (2,\ymax)  -- (4,\ymax) ;
\draw[thick]        (4,\ymax)  -- (12,\ymax) ;
\draw[thick,dashed] (12,\ymax) -- (14,\ymax) ;

\draw[thick] (5,\ymin)  -- (5,\ymax)  ;
\draw[thick] (7,\ymin)  -- (7,\ymax)  ;
\draw[thick] (9,\ymin)  -- (9,\ymax)  ;
\draw[thick] (11,\ymin) -- (11,\ymax) ;

\def\xcenter{8}
\draw[thick,->] (\xcenter,6.5) -- (\xcenter,5.5) ;
\def\xcenter{6}
\draw[thick,->] (\xcenter,6.5) -- (5,5.5) ;
\def\xcenter{10}
\draw[thick,->] (\xcenter,6.5) -- (11,5.5) ;

\def\ymin{1}
\FPadd\ymax{\ymin}{\Height}

\draw[thick] (4,\ymin) -- (6,\ymin) ;
\draw[thick] (4,\ymax) -- (6,\ymax) ;
\def\xloc{3.75}
\fill[gray!40] (\xloc,\ymin) rectangle (4,\ymax) ;
\def\xloc{6.25}
\fill[gray!40] (\xloc,\ymin) rectangle (6,\ymax) ;

\draw (5,\ymin) node[above] { $P_{k-1}$ } ;
\draw[thick] (7,\ymin) -- (9,\ymin) ;
\draw[thick] (7,\ymax) -- (9,\ymax) ;
\def\xloc{6.75}
\fill[gray!40] (\xloc,\ymin) rectangle (7,\ymax) ;
\def\xloc{9.25}
\fill[gray!40] (9,\ymin)     rectangle (\xloc,\ymax) ;
\draw (8,\ymin) node[above] { $P_{k}$ } ;

\draw[thick] (10,\ymin) -- (12,\ymin) ;
\draw[thick] (10,\ymax) -- (12,\ymax) ;
\def\xloc{9.75}
\fill[gray!40] (\xloc,\ymin)  rectangle (10,\ymax) ;
\def\xloc{12.25}
\fill[gray!40] (12,\ymin) rectangle (\xloc,\ymax) ;
\draw (11,\ymin) node[above] { $P_{k+1}$ } ;

\end{tikzpicture}
  \end{center}
  \caption{Splitting of the computational grid among the MPI processes $P_k,\, k=0,\ldots,p-1$. The gray zone refers to additional (ghost points) storage used for the MPI communications
    between neighbooring processes.}
  \label{grid_partition}
\end{figure}
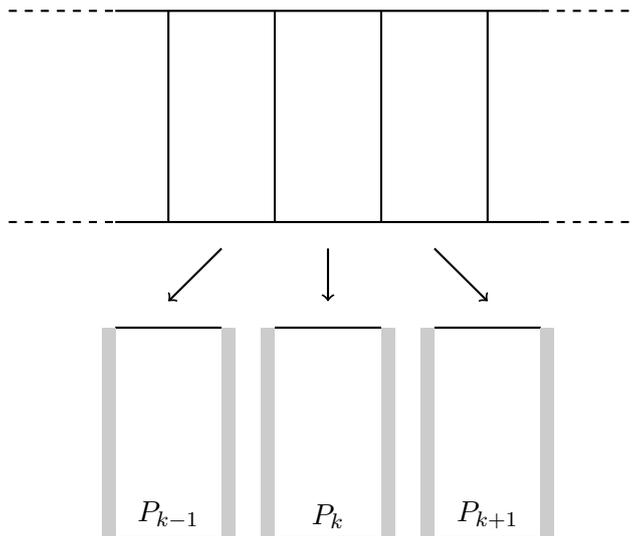

The first step when developping a parallel algorithm is to define a suitable and efficient splitting of the datas among the MPI processes~:
each MPI process will treat datas associated with a part of the total computational mesh. For our problem, this choice is straightforward and is related to the algorithm
used to solve the linear systems. Indeed, it is much easier to implement a parallel resolution of tridiagonal linear systems rather than a parallel version
of the DFT. Therefore, the parallel version of the code is based on a splitting of the datas along the horizontal axis, so that each MPI process
works with a vertical slice of the computational mesh as it is illustrated on Figure~\ref{grid_partition}.
%
%
%
%

In the framework of finite volume or finite difference schemes on cartesian grids, the explicit computation of spatial derivatives is local and involves
very few communications. The only tricky part concerns the resolution of the linear systems.
The step ii) of the direct solver described in the previous section consists in solving a linear system involving the matrix $I_{n_c}\,+\, M\,G^{-1}Q$.
As the $LU$-factorization of this matrix has been computed and stored in a pre-processing step at the beginning of the time iterations, we have to solve
two triangular systems which can not be efficiently performed on parallel computers. As $n_c$ is small compared to the size of the global problem, we
choose to dedicate this task to one given MPI process (fixed in advance) per unknown, that is $u,v$ and $p.$ Once these linear systems are solved,
the resulting vectors are scattered from these MPI processes to the other ones.

As it was previously mentioned, linear systems of steps i) and ii) are solved by first applying a DFT in the $y$-direction~: these computations
are independant and can be performed without any communications due to the distribution of datas among the MPI processes. This results in a collection, one per
grid point in the $y$-direction, of independant tridiagonal linear systems connecting all nodes of the mesh in the $x$-direction.
A parallel direct solver based on the \textit{divide and conquer approach} (DAC) for tridiagonal matrices has been implemented (see~\cite{Bondeli1991}).
The DAC method, applied to solve one tridiagonal linear system on $n_p>1$ MPI processes, consists in splitting the tridiagonal matrix into $n_p$ independant blocks (one per MPI process).
The solutions of these systems have to be corrected in order to recover the solution of the global system. These corrections correspond to $2n_p-1$ values which are solutions of a tridiagonal linear system
of size $2n_p-1$. This phase of the DAC method is sequential and has to be performed on one process inducing a useless waiting time for the other processes.
However, as we have to solve $n_y$ such systems simultaneously, this sequential part can be distributed among all the $n_p$ processes. In this context, the DAC algorithm
leads to an efficient parallel code.

The parallel code has a good level of performance~: less than $15\%$ of the CPU time is spent in communications between MPI processes. The sequential
part performed on one process represents a negligible amount of CPU time. The computations presented here have been performed on a DELL cluster using up to 32 cores of
Xeon processors. A low latency bandwith network connects the cluster nodes.

\subsection{Iterative solver for the case of moving bodies}
For solid bodies moving in a computational domain filled by a fluid, as the case considered in Section~\ref{Num_res_moving_bodies}, cut-cells may change at each time iteration. Therefore,
the coefficients of the matrices of the linear systems for the velocity components, issued from the discretizetion of the momentum equations, and for the pressure increment computed in the
projection step of the time scheme, have to be recomputed at each time step. In that context, the use of the fast direct solver described in the previous section is cumbersome and inefficient
except on coarse meshes. In order to be able to treat such configurations, we have implemented a PETSc version~\cite{petsc-web-page,petsc-user-ref} of our cut-cell scheme.
The main advantage of the PETSc Library is that, in a parallel programming environment based on MPI, many iterative solvers combined with different preconditionners can be used.
The choice can be made at run time.

\section{Numerical Results}
\label{sec:numer}

\subsection{Flow past a circular cylinder at $\textrm{Re}=1000$ and $3000$}
In this section, we present numerical simulations, performed with the parallel direct solver method described in Section 3.1. We consider the case of flows past a fixed circular cylinder of 
diameter $D$. The Reynolds number is defined based on the diameter $D$ of the cylinder, \textit{i.e.} $\textrm{Re}=U_{\infty}D/\nu$ where $U_{\infty}$ is the horizontal free stream velocity.
As non-dimensional time, we consider $T=2U_{\infty}t/D$. 

A circular cylinder of diameter equals to unity is centered at the origin of the computational domain $\Omega=(-L_x,L_x)\times (-L_y,L_y)$.
As boundary conditions, a uniform velocity profile $\mathbf{u}(\mathbf{x}=-L_x,t)=(1,0)$ is imposed at the inflow and a convective boundary condition is applied at the exit,
namely the convective equation
\begin{equation}
  \frac{\partial \mathbf{u}}{\partial t} + (1,0)\cdot \bm{\nabla}\mathbf{u} = 0
\end{equation}
is solved at $x=L_x$. On the top and bottom boundaries, that is $y=\pm L_y$, slip boundary conditions are used that is $\frac{\partial u}{\partial n}=0$ and $v=0$. Finally,
no-slip ($\mathbf{u}|_{\Gamma_S}=\mathbf{0}$) boundary condition is applied on the surface of the obstacle. 

\begin{figure*}
\begin{center}
  \includegraphics[width=0.8\textwidth]{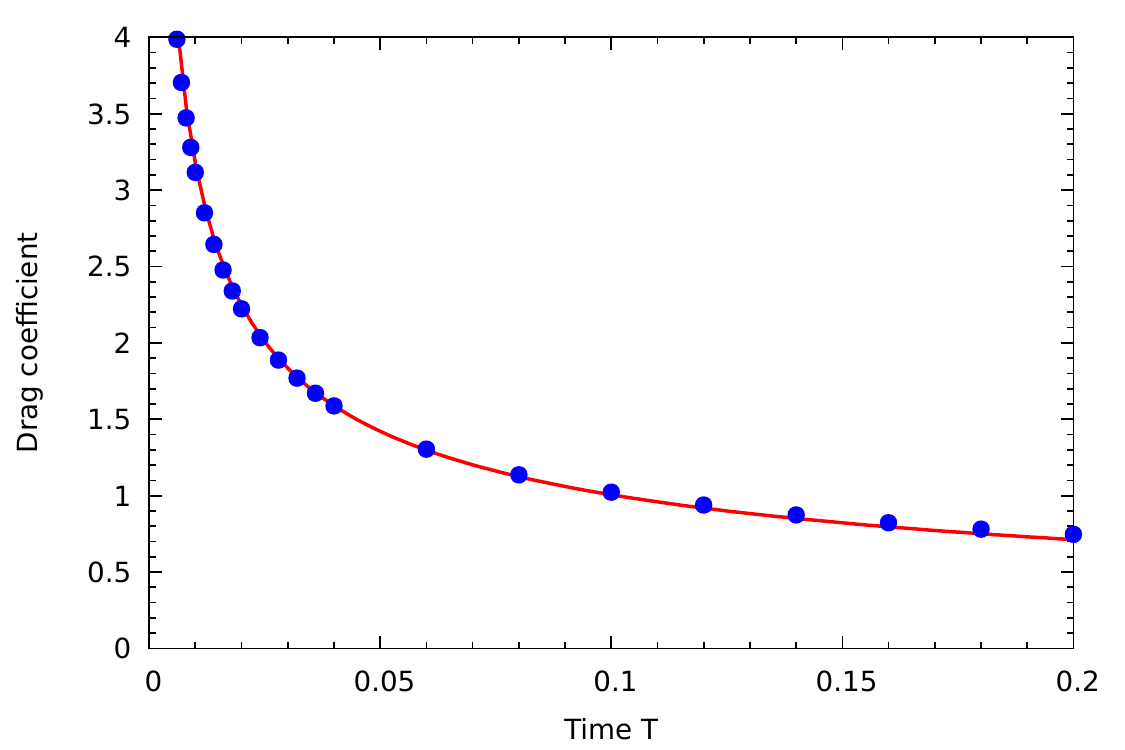}
\end{center}
\caption{Evolution of the drag coefficient of a circular cylinder at $\textrm{Re}=1000$. Solid line (red): theoretical prediction~\eqref{Theoric_drag}; Blue dots~: numerical results on a $4096\times 8192$ grid in $\Omega=(-10,10)^2$.}
\label{fig:Drag_sqrs_1000}
\end{figure*}

\begin{figure*}
\begin{center}
  \includegraphics[width=0.8\textwidth]{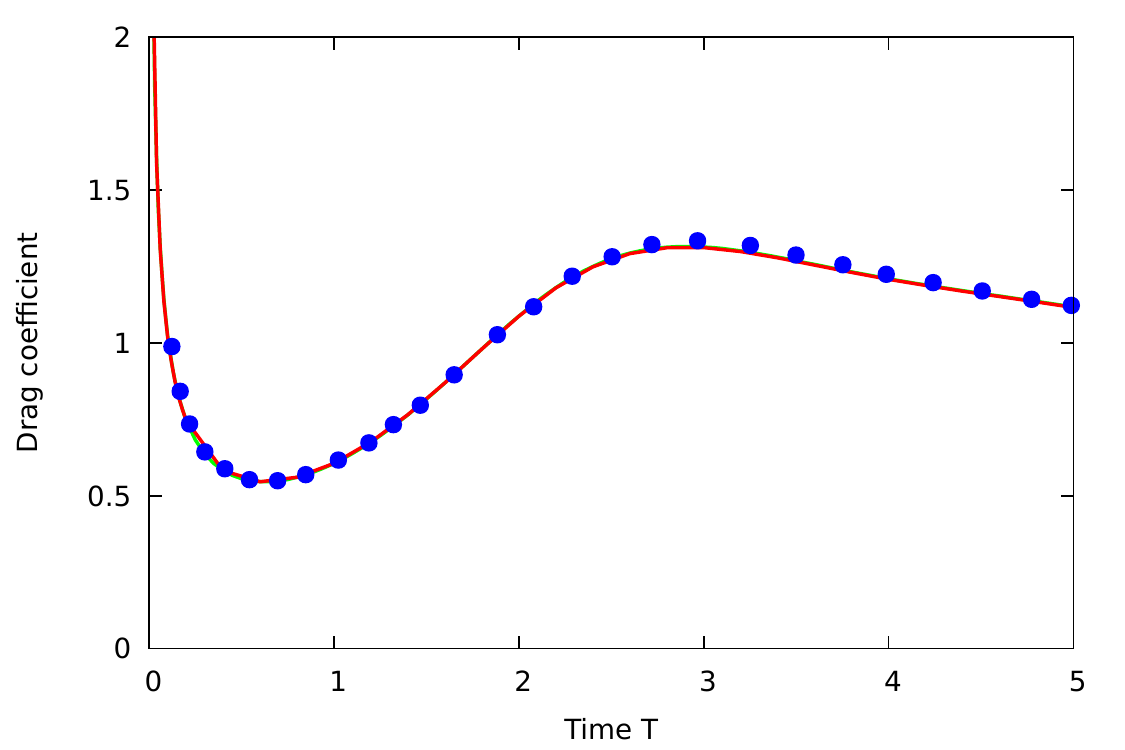}
\end{center}
\caption{Evolution of the drag coefficient of a circular cylinder at $\textrm{Re}=1000$. Red solid line: numerical results on a $4096\times 8192$ grid in $\Omega=(-10,10)^2$; Green solid line: numerical results on a $8192\times 16384$ grid in $\Omega=(-10,10)^2$; Blue dots: results from~\cite{KL1995}.}
\label{fig:Drag_1000}
\end{figure*}

\begin{figure*}
  \begin{center}
    \begin{tabular}{cc}
      \includegraphics[width=0.45\textwidth]{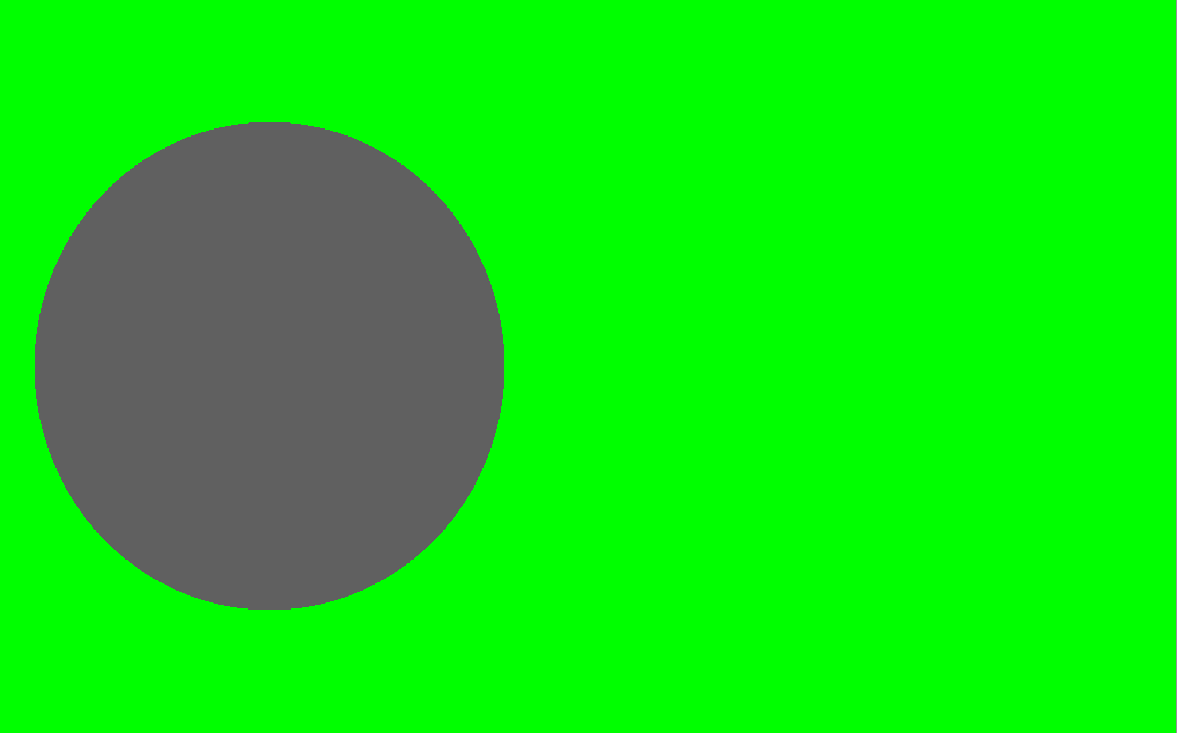} & \includegraphics[width=0.45\textwidth]{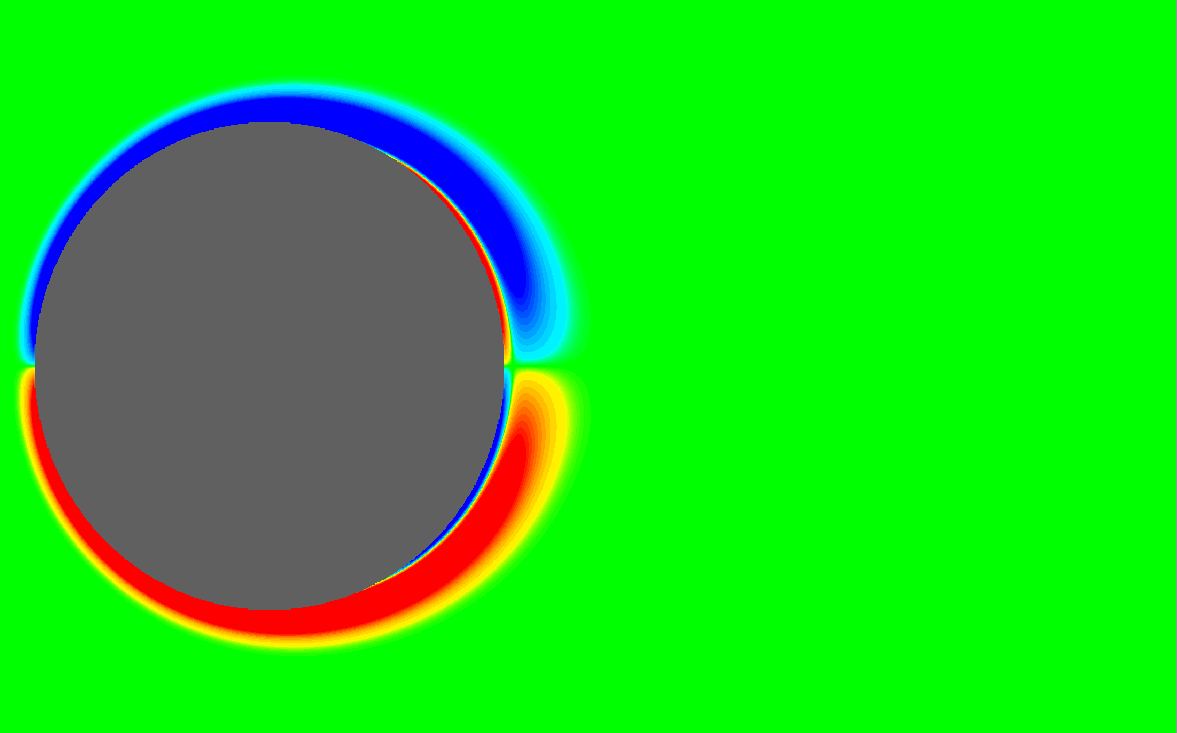} \\
      $T=0$ &  $T=1$ \\
      \includegraphics[width=0.45\textwidth]{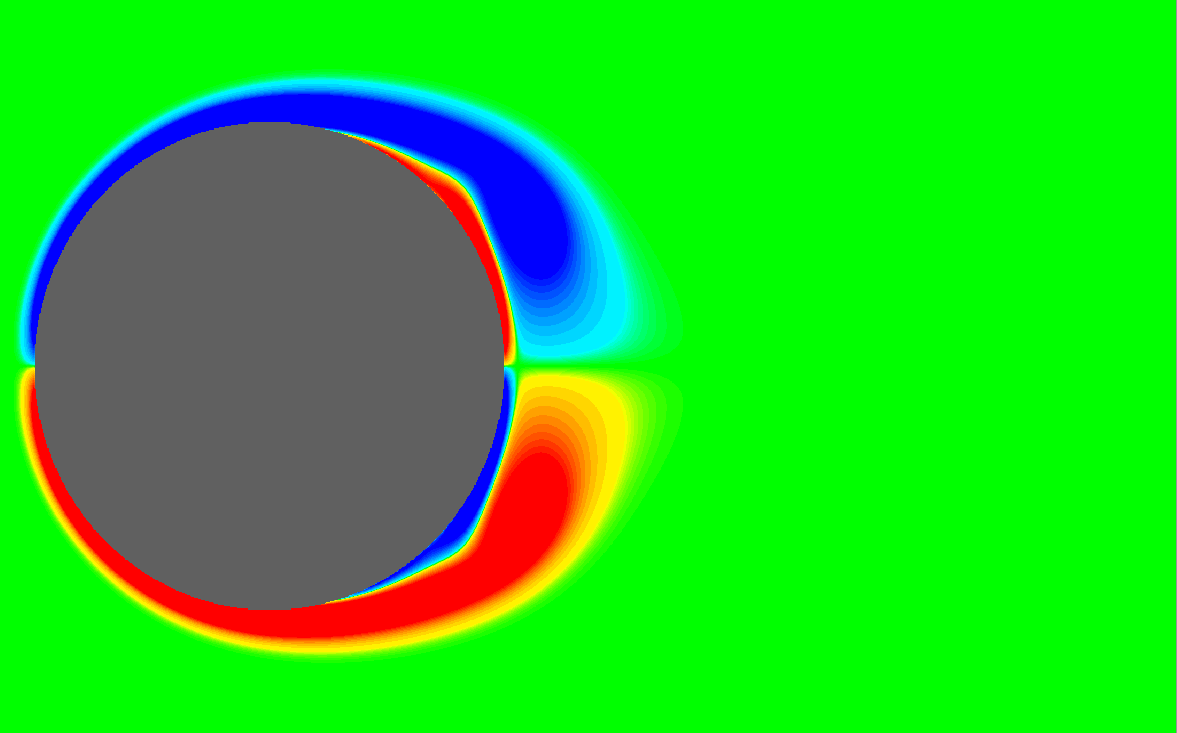} &  \includegraphics[width=0.45\textwidth]{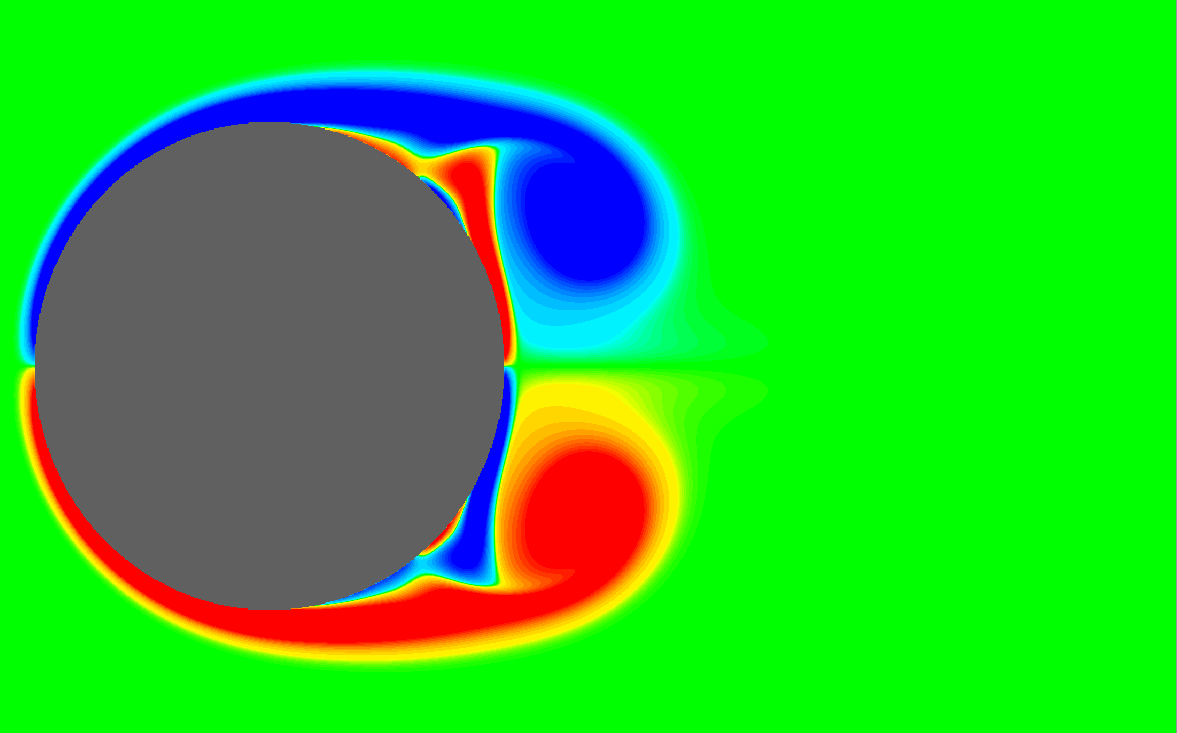} \\
      $T=2$ &  $T=3$ \\
      \includegraphics[width=0.45\textwidth]{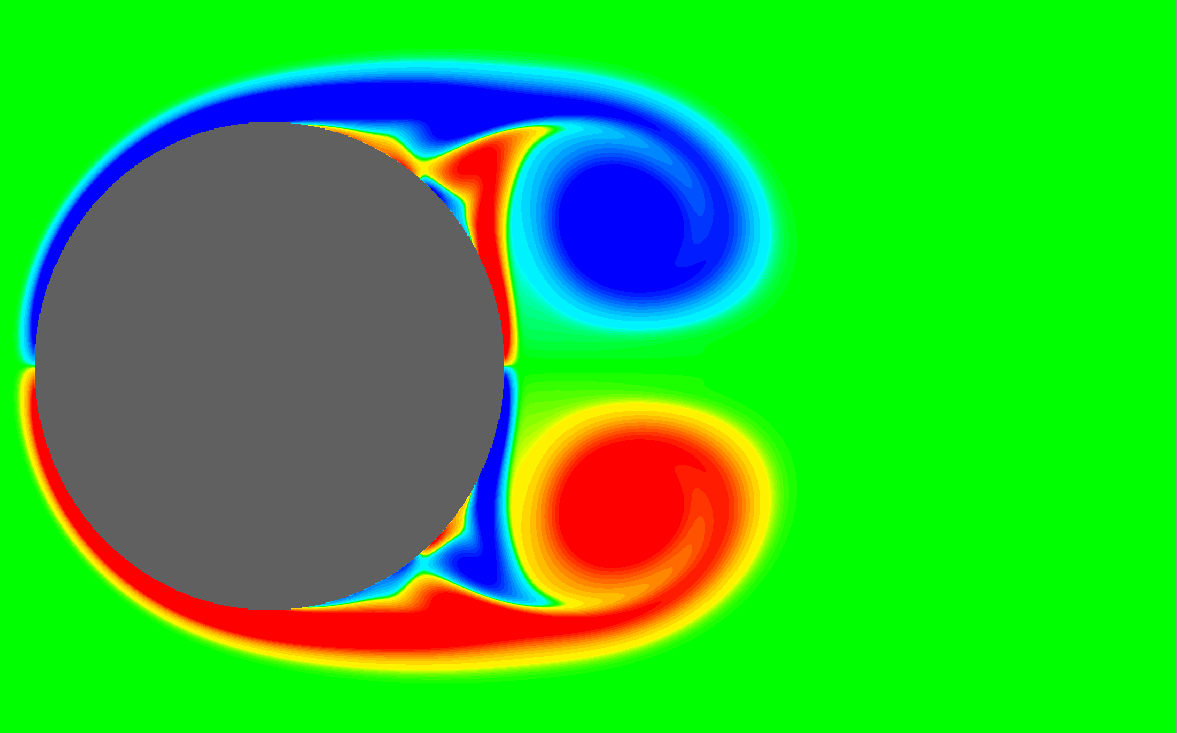} &  \includegraphics[width=0.45\textwidth]{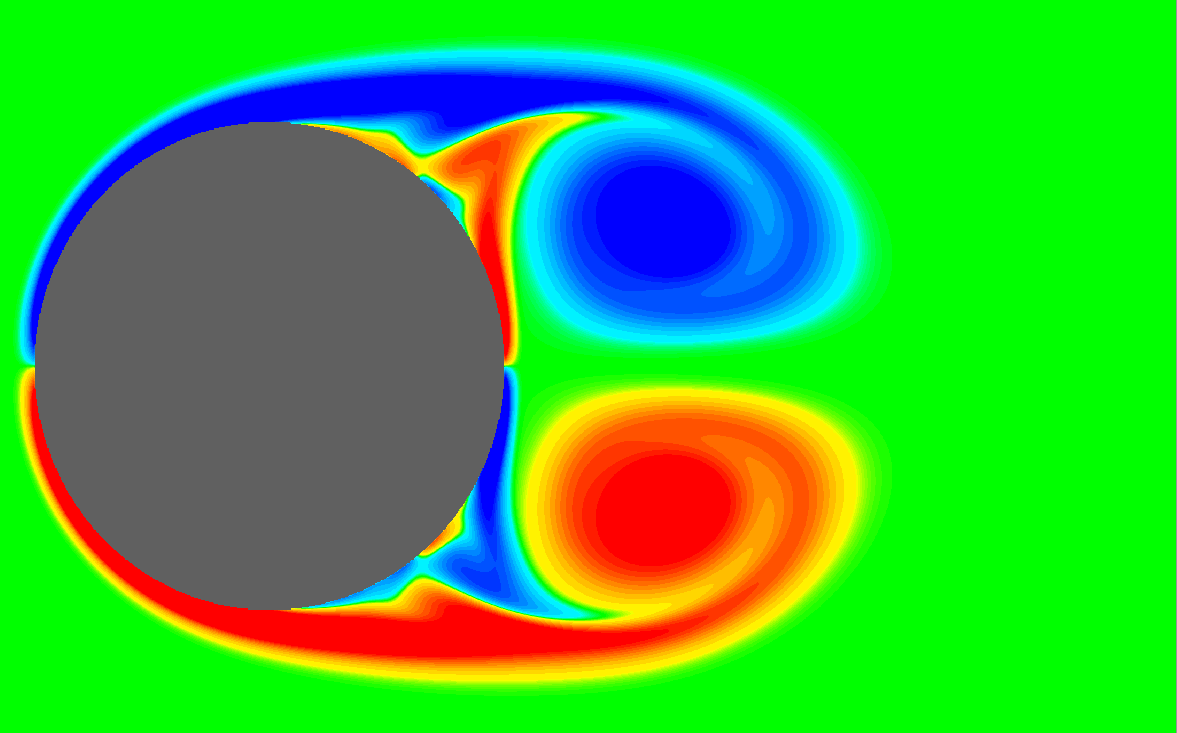} \\
      $T=4$ &  $T=5$ \\
      \includegraphics[width=0.45\textwidth]{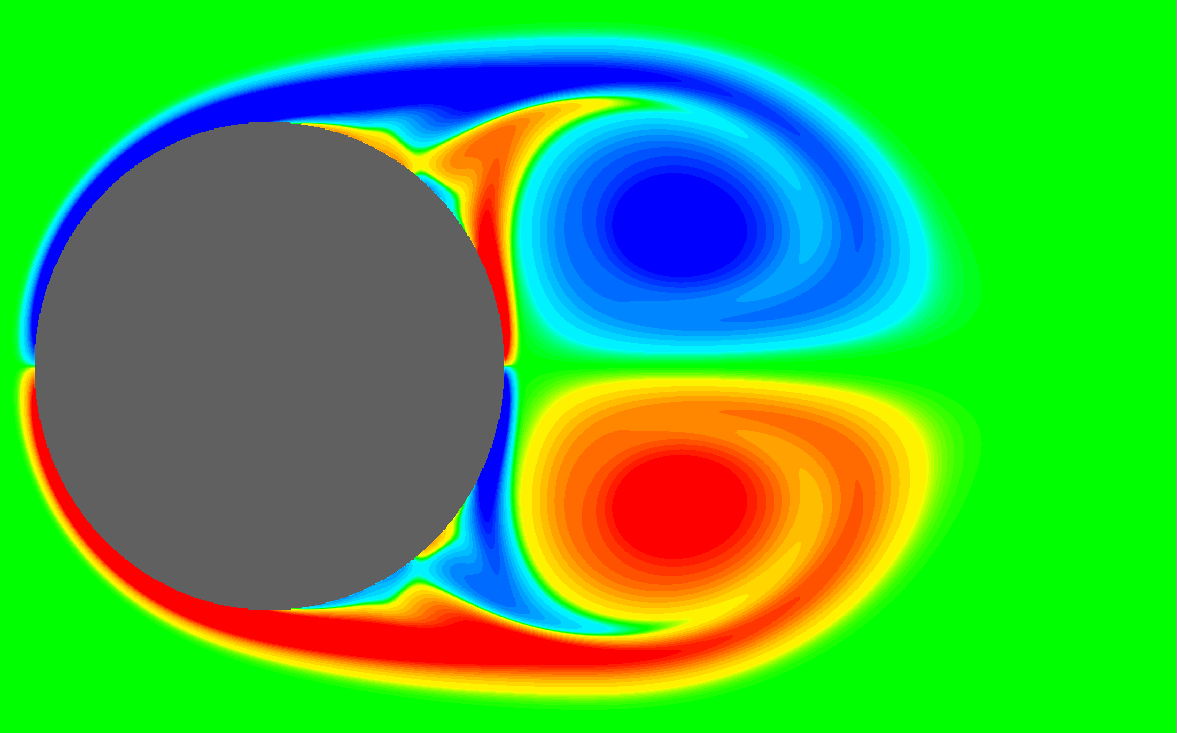} &  \includegraphics[width=0.45\textwidth]{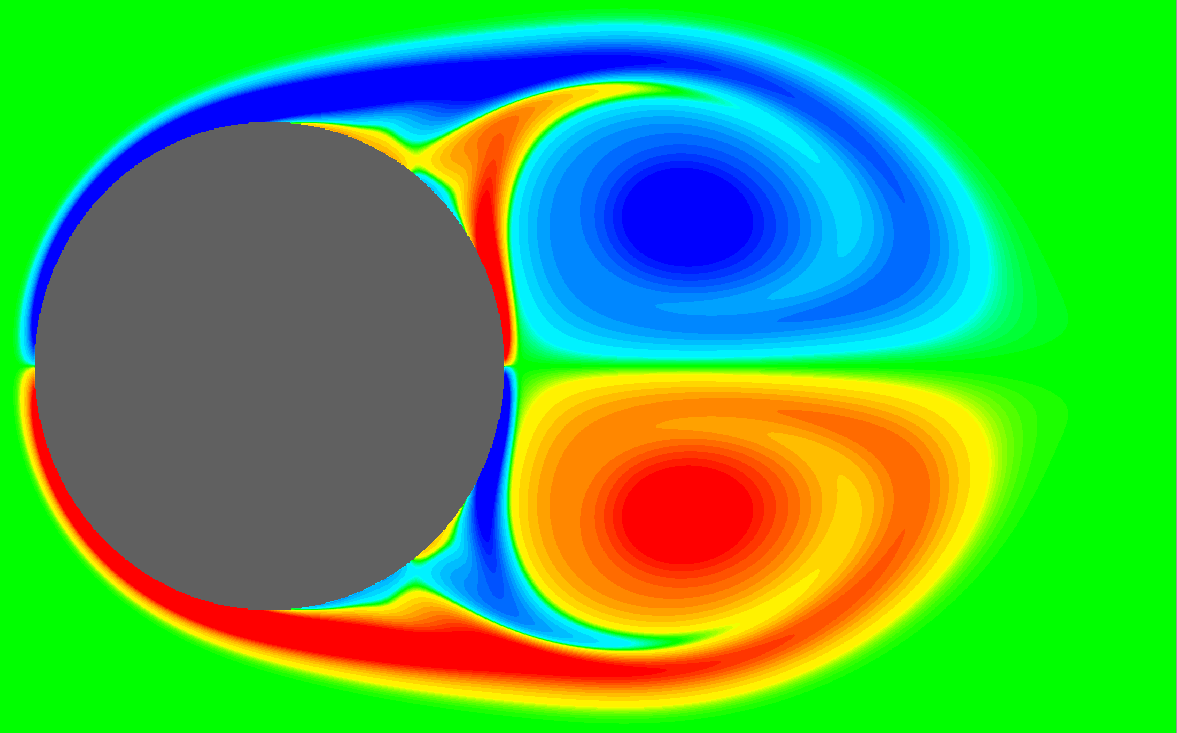} \\
      $T=6$ &  $T=7$ \\
    \end{tabular}
\end{center}
\caption{Vorticity of a flow past a circular cylinder at $\textrm{Re}=1000$ simulated on a grid with $4096\times 8192$ mesh points in the computational domain $\Omega=(-10,10)^2$ at different times $T\in[0,7]$.}
\label{fig:Cylinder_1000}
\end{figure*}

For this problem important quantities reflecting the dynamics of the vortices formed in the vicinity of the solid boundary and developping at the rear of the cylinder are the pressure drag and lift coefficients.
They are derived from the total drag force on the body, which is computed as
\begin{equation}
  \mathbf{F}_b\,=\, \int_{\Gamma_S} \Bigl(-p\mathbf{n}+\frac 1{\textrm{Re}}\frac{\partial\mathbf{u}}{\partial n} \Bigr) ds.
\end{equation}
The pressure drag and lift coefficients $C_p$ and $C_\ell$ are given by $C_p=2\mathbf{F}_b\cdot\mathbf{e}_x$ and $C_\ell=2\mathbf{F}_b\cdot\mathbf{e}_y$ and the (total) drag coefficient is $C_d=C_p+C_\ell$.
Starting with a flow at rest, the drag coefficient behaves as $T^{-1/2}$ in the early stage of the development of the vortices. This square-root singularity has been theoretically predicted
by Bar-Lev and Yang in~\cite{BY}. They have derived the following expression for the (total) drag coefficient
\begin{equation} \label{Theoric_drag}
  C_{\textrm{pred}} = 4\sqrt{\frac{2\pi}{\textrm{Re}T}} + \frac{2\pi}{\textrm{Re}}\left( 9 -\frac{15}{\sqrt{\pi}}\right).
\end{equation}
As whown in Figure~\ref{fig:Drag_sqrs_1000}, the values obtained with the cut-cell scheme on a grid with $4096\times 8192$ mesh points discretizing the domain $\Omega=(-10,10)^2$
perfectly match the theoretical curve drawing~\eqref{Theoric_drag} on the time interval $T\in[0,0.2]$ for $\textrm{Re}=1000$. Our cut-cell method captures the square-root singularity of
the drag coefficient. This simulation has been run using $16$ MPI processes. On longer time interval, namely $T\in [0,5]$, the results are in good agreement with those obtained by Koumoutsakos
and Leonard in~\cite{KL1995} with a vortex method. In order to test the grid convergence of these results, the same simulation has been conducted on a grid with two times more points in
both spatial directions, that is $8192\times 16384$ mesh points, in the same computational domain. In that case, $32$ MPI processes have been used. Both results are almost indistinguishable on
Figure~\ref{fig:Drag_1000} indicating that the coarser resolution is enough to capture the essential features of the flow at this Reynolds number.
The development of the flow around the impulsively started cylinder at $\textrm{Re}=1000$ can be seen on Figure~\ref{fig:Cylinder_1000}.
In the early stage $T\le 1$, a primary vortex develops in the vinicity of the boundary at the rear of the cylinder.
Then for $T\in[1,2]$ a secondary vortex appears trying to move insight the primary vortex and to push it away from the solid boundary ($T\ge 4$). 
A tertiary vortex is visible at $T=3$ which remains sticked to the boundary constrained by the two other vortices having more strength.
These results compare well with the same flow representations shown in~\cite{KL1995}. As expected on short time interval ($T\le 5$) the flow remains symmetric. 
\begin{figure*}
  \begin{center}
  \includegraphics[width=0.8\textwidth]{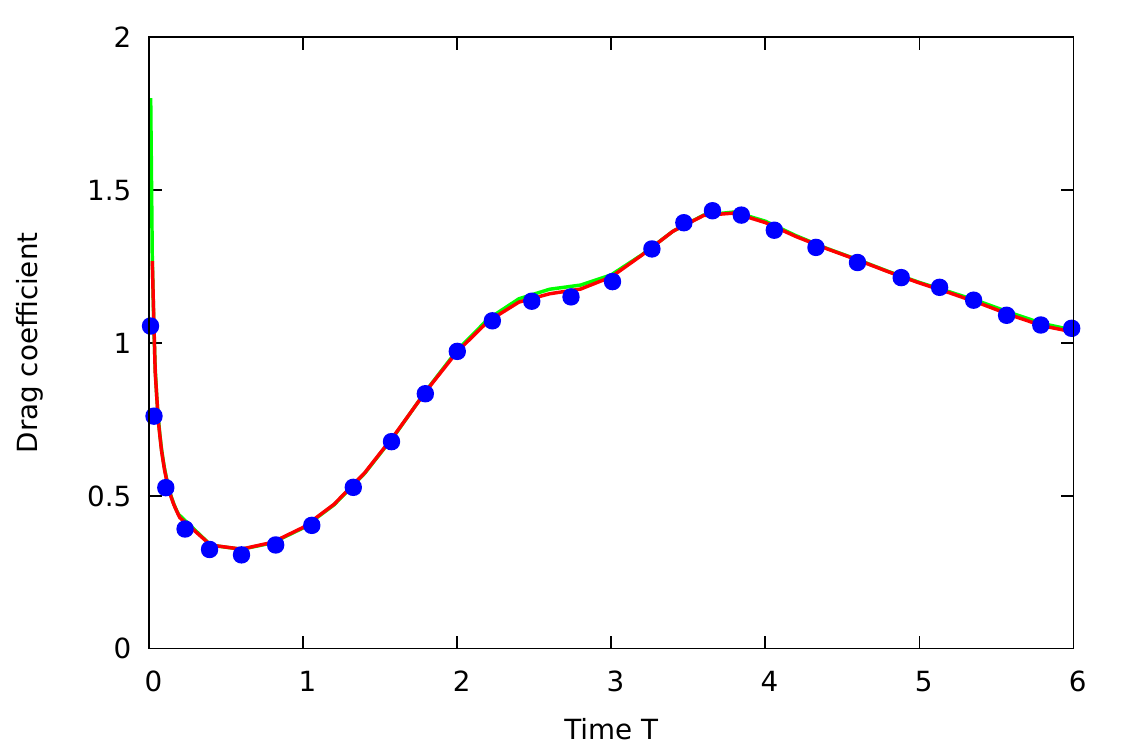}
\end{center}
\caption{Evolution of the drag coefficient of a circular cylinder at $\textrm{Re}=3000$. Red solid line: numerical results on a $4096\times 8192$ grid in $\Omega=(-10,10)^2$; Green solid line: numerical results on a $8192\times 16384$ grid in $\Omega=(-10,10)^2$; Blue dots: results from~\cite{KL1995}.}
\label{fig:Drag_3000}
\end{figure*}

\begin{figure*}
  \begin{center}
    \begin{tabular}{cc}
      \includegraphics[width=0.45\textwidth]{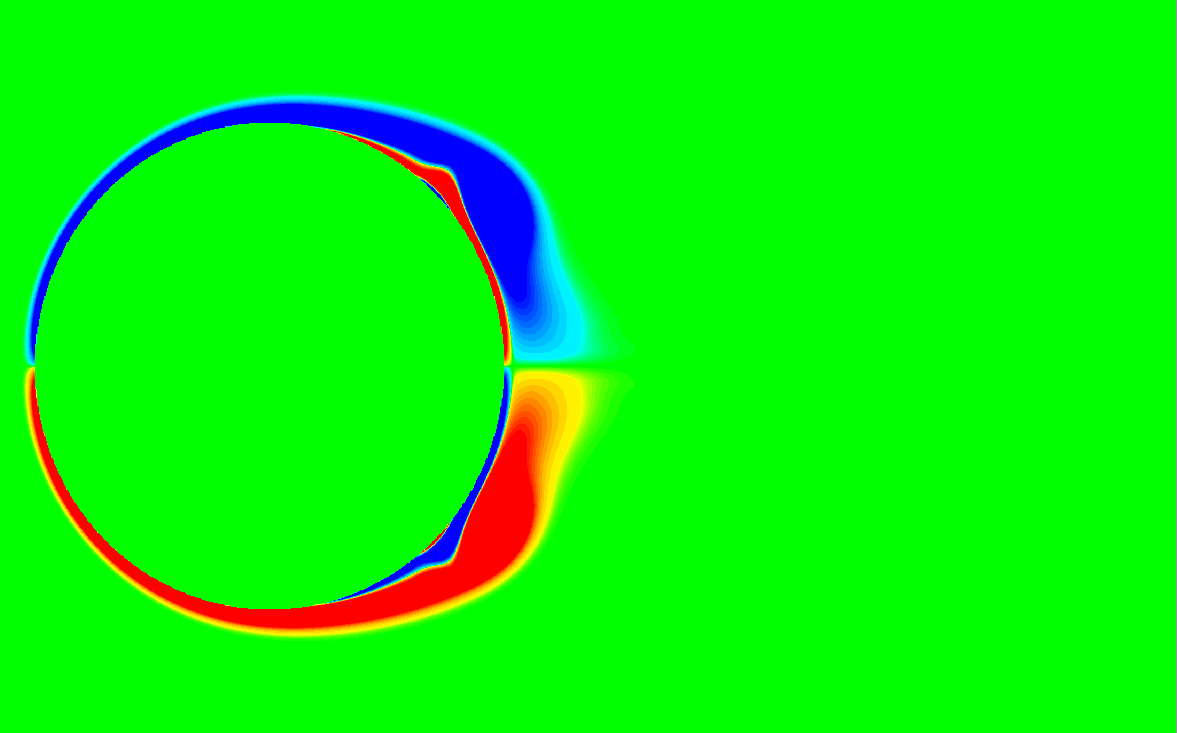} &  \includegraphics[width=0.45\textwidth]{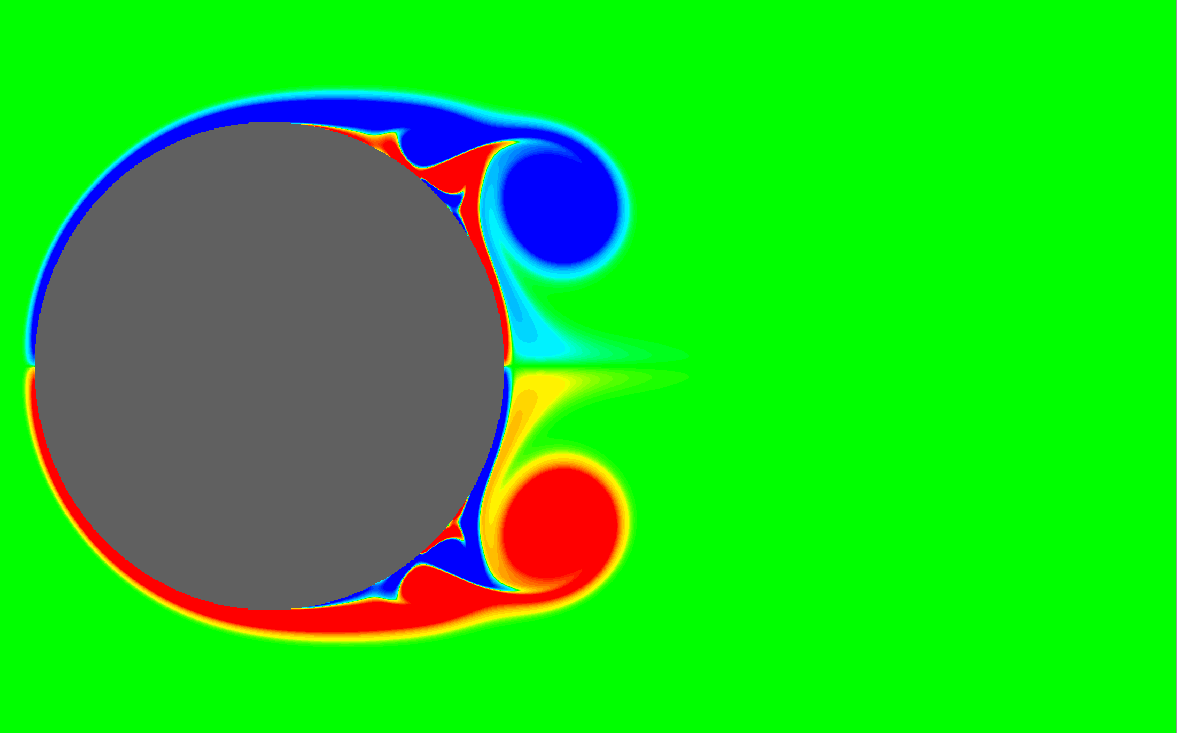} \\
      $T=2$ &  $T=3$ \\
      \includegraphics[width=0.45\textwidth]{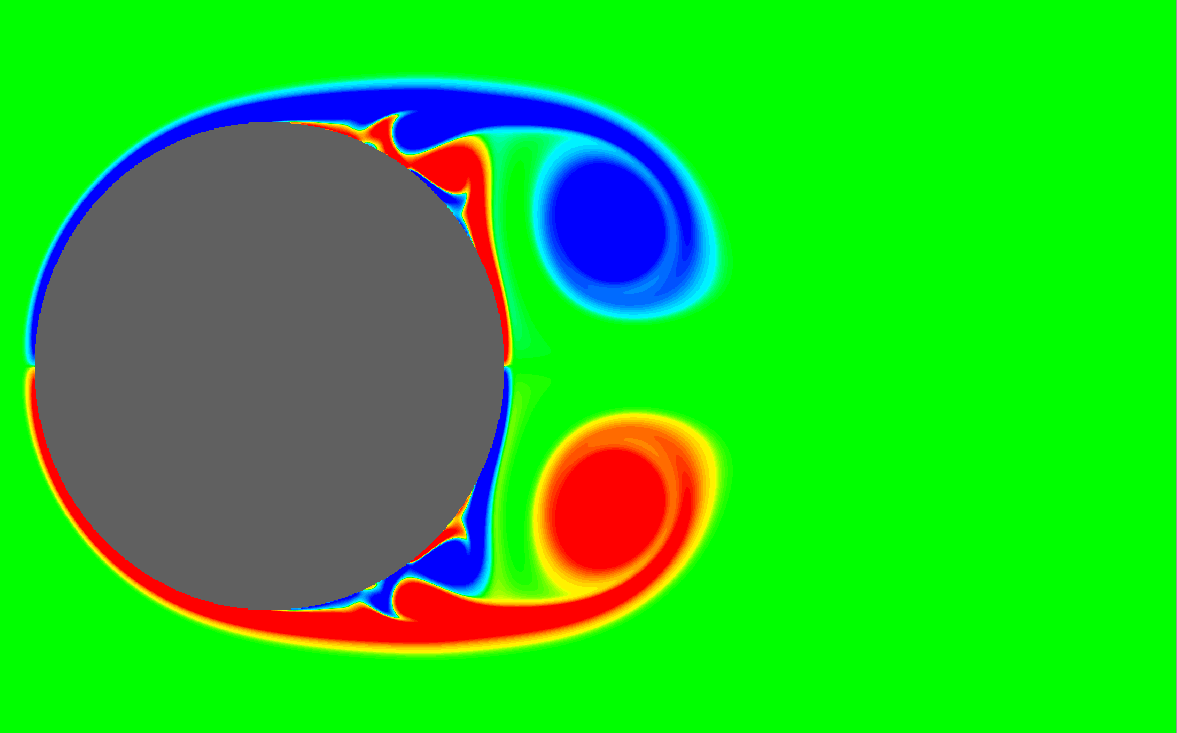} &  \includegraphics[width=0.45\textwidth]{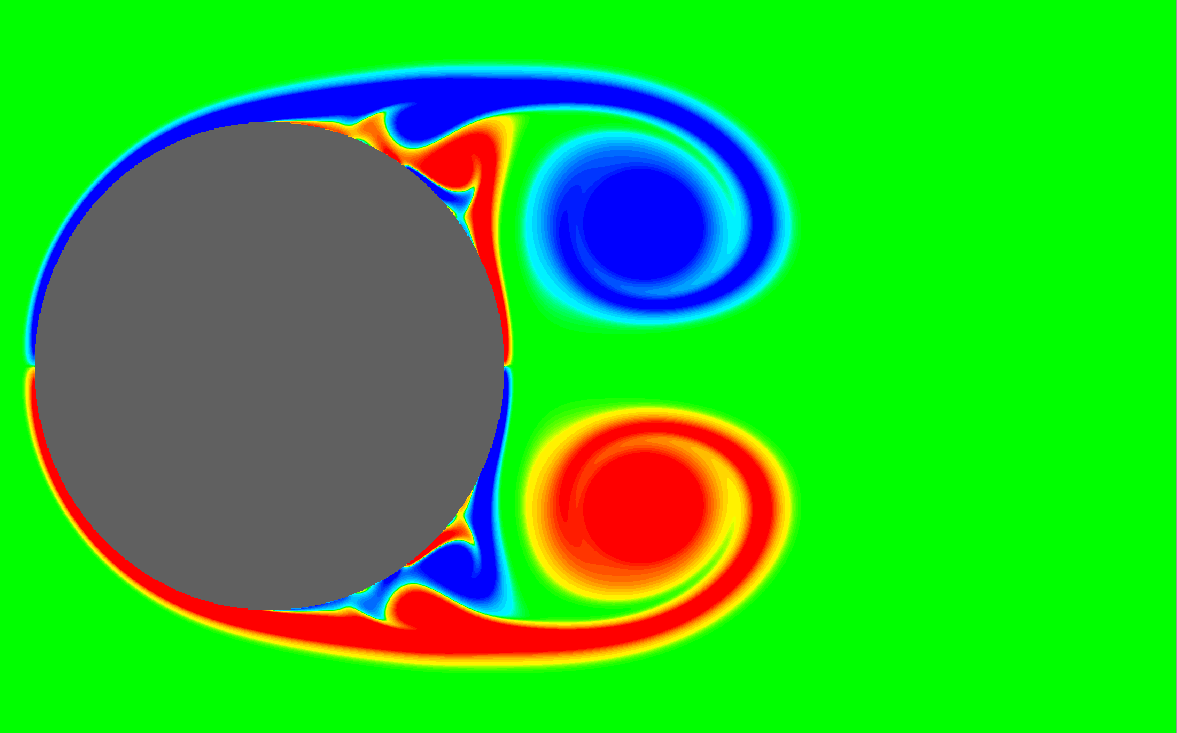} \\
      $T=4$ &  $T=5$ \\
      \includegraphics[width=0.45\textwidth]{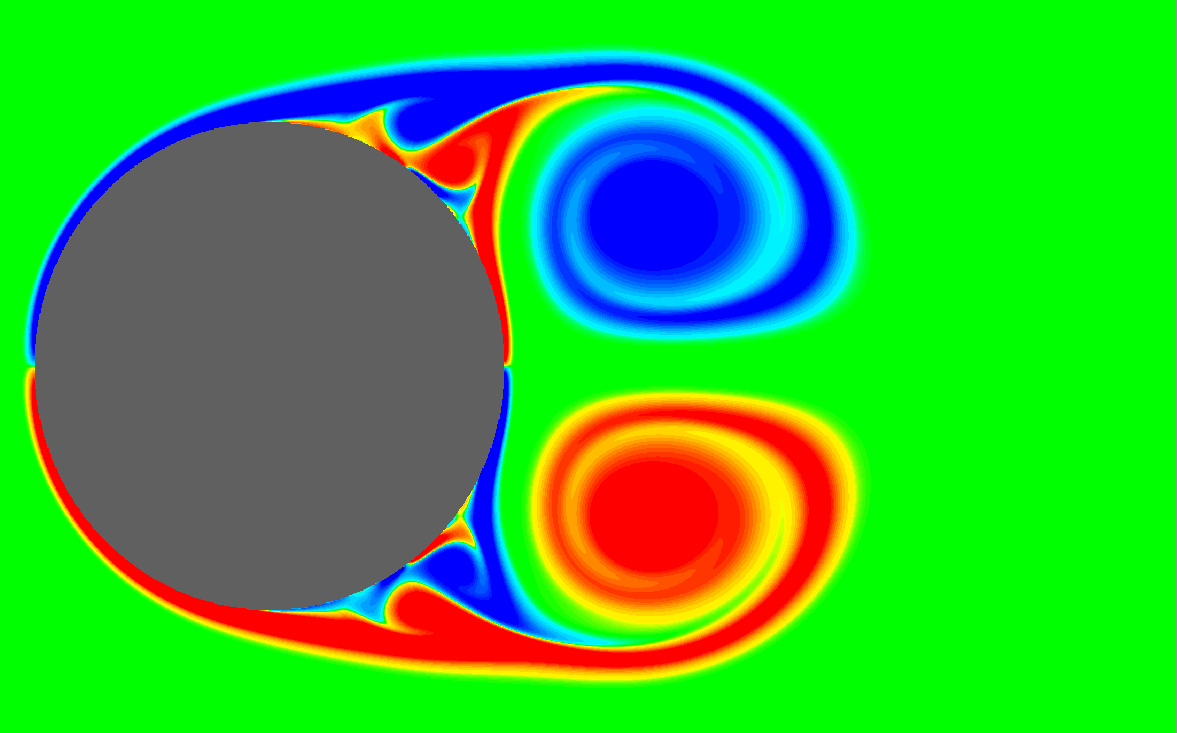} &  \includegraphics[width=0.45\textwidth]{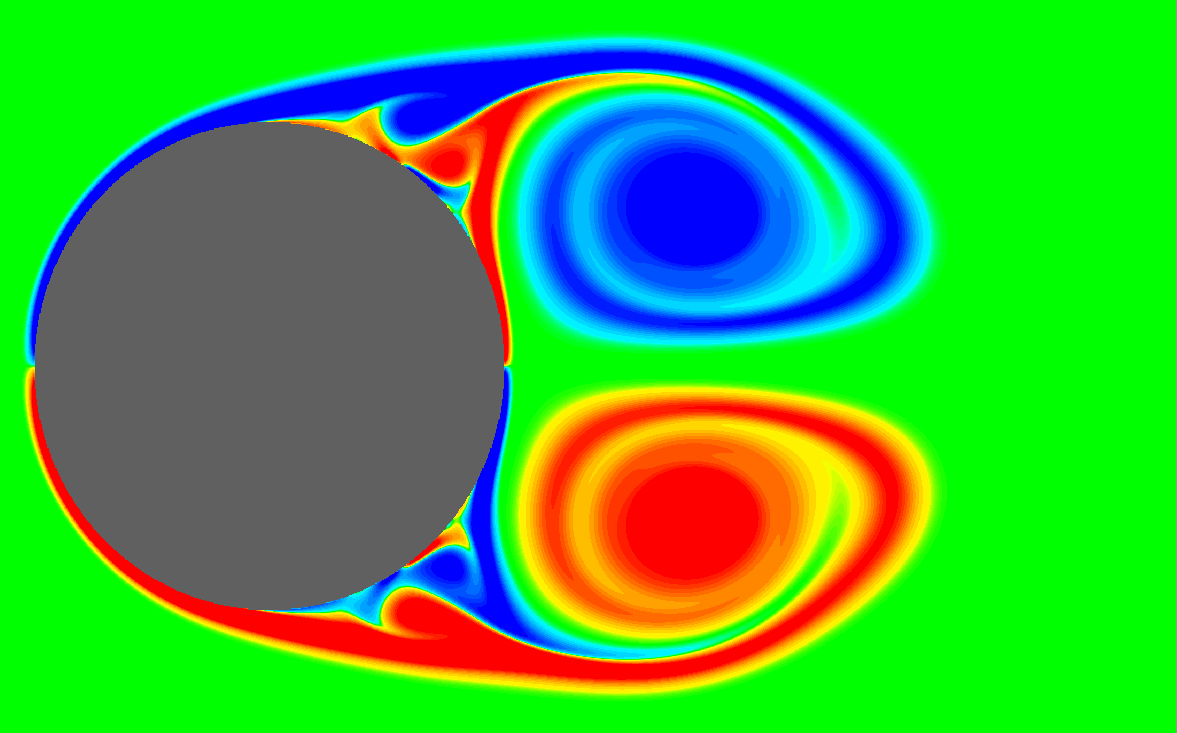} \\
      $T=6$ &  $T=7$ \\
      \includegraphics[width=0.45\textwidth]{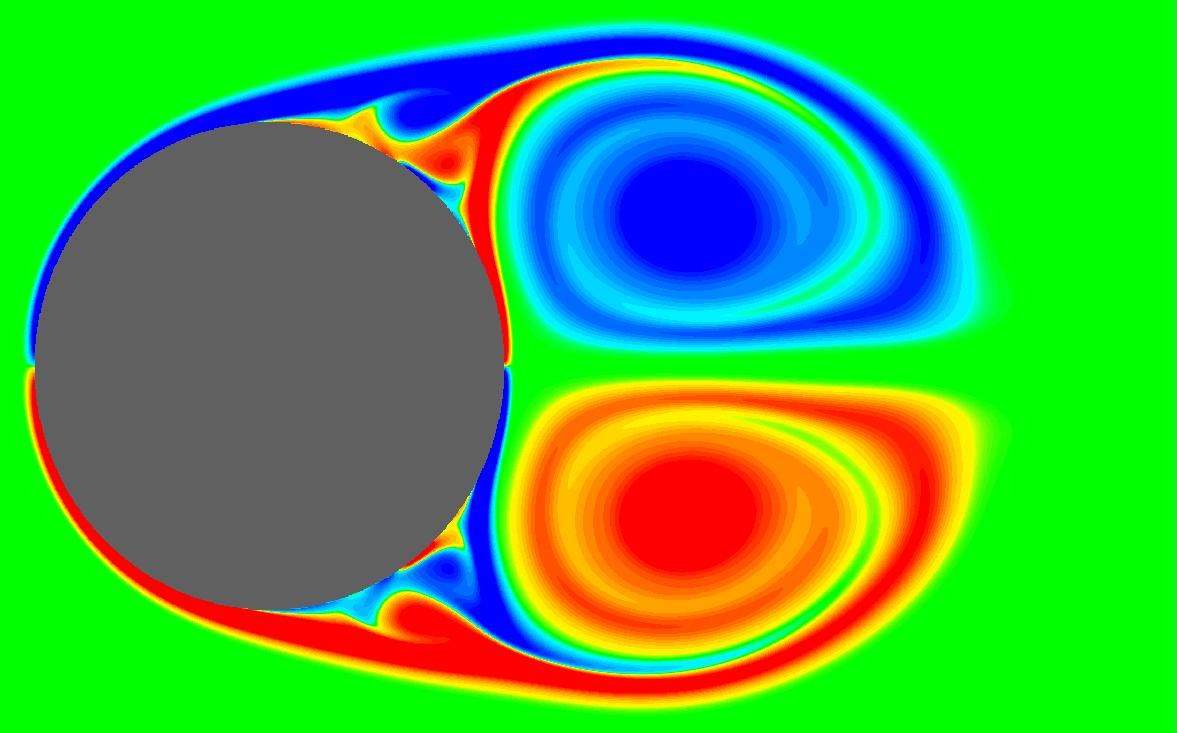} &  \includegraphics[width=0.45\textwidth]{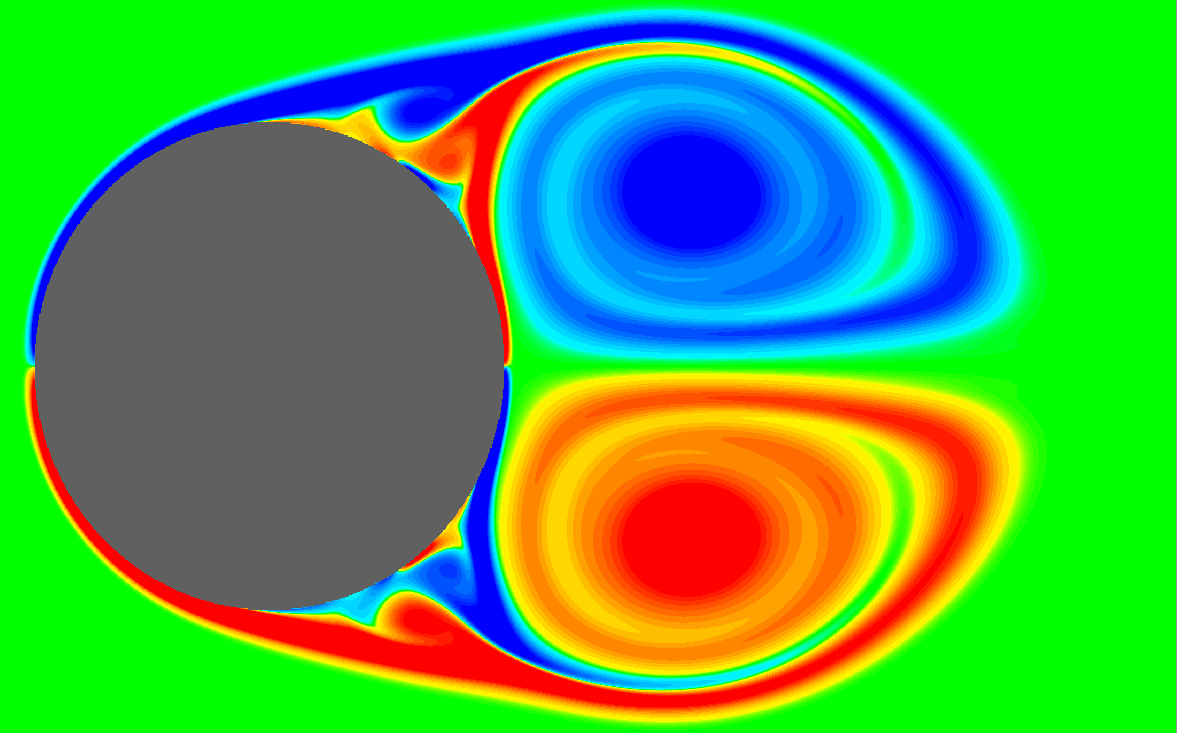} \\
      $T=8$ &  $T=9$ \\
    \end{tabular}
\end{center}
\caption{Flow around a circular cylinder at $\textrm{Re}=3000$ simulated on a grid with $4096\times 8192$ mesh points in the computational domain $\Omega=(-10,10)^2$ at different times $T\in[0,9]$.}
\label{fig:Cylinder_3000}
\end{figure*}

At $\textrm{Re}=3000$, the time evolution of the drag coefficient plotted on Figure~\ref{fig:Drag_3000} exhibits also the square-root singularity on short time interval
and is in good agreement with the results of Koumoutsakos and Leonard. As expected, the drag coefficients remains almost constant for $T\in[2,3]$ (see~\cite{KL1995}).
Note that on this time interval, a small difference exists between the results computed on the two different grids. However, the coarser simulation is fine enough to
capture the flow dynamics at $\textrm{Re}=3000$. The mesh size of the coarser grid, which is constant in the vicinity of the solid boundary, is $h=20/8192\approx 2.44\times 10^{-3}$.
Note that the size of the computational domain and the boundary conditions imposed at the exit may influence the results. Estimating the values for $L_x$ and $L_y$ required so that the numerical
results being of the order of the numerical scheme error, namely $\textrm{O}(h^2)$, is an open question. This will be addressed in further works.
At this Reynolds number, a scenario similar than that at $\textrm{Re}=1000$ can be observed on Figure~\ref{fig:Cylinder_3000} with the development of three vortices in the early
stage of the flow dynamics. The secondary vortex penetrates further inside the primary vortex aera and the tertiary vortex has more strength as it could be expected with
less effects of the viscous forces at $\textrm{Re}=3000$. Again, an overall good agreement is found with the vortcity contours shown in~\cite{KL1995} at the same Reynolds number and time $T$.

\begin{figure*}
\begin{center}
  \includegraphics[width=0.8\textwidth]{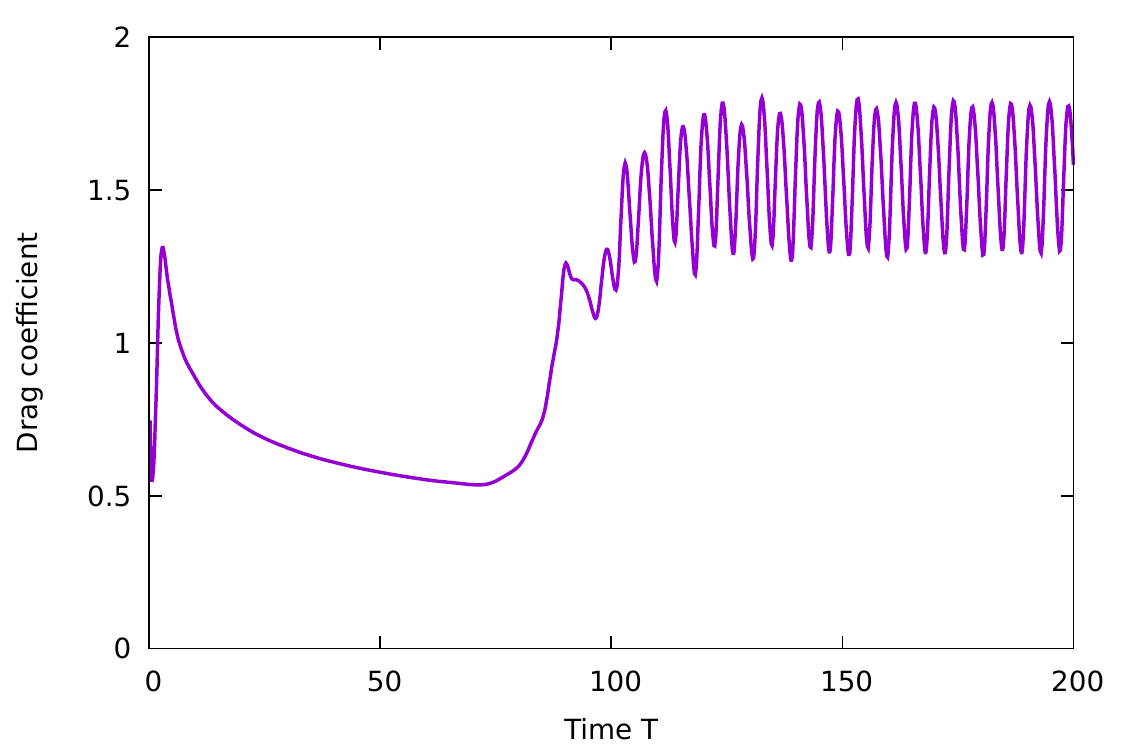}
\end{center}
\caption{Evolution of the drag coefficient of a circular cylinder for $\textrm{Re}=1000$ simulated on a grid with $4096\times 8192$ mesh points in the computational domain $\Omega=(-10,10)^2$. }
\label{fig:Drag_1000_Time_Long}
\end{figure*}

\begin{figure*}
\begin{center}
  \includegraphics[width=0.8\textwidth]{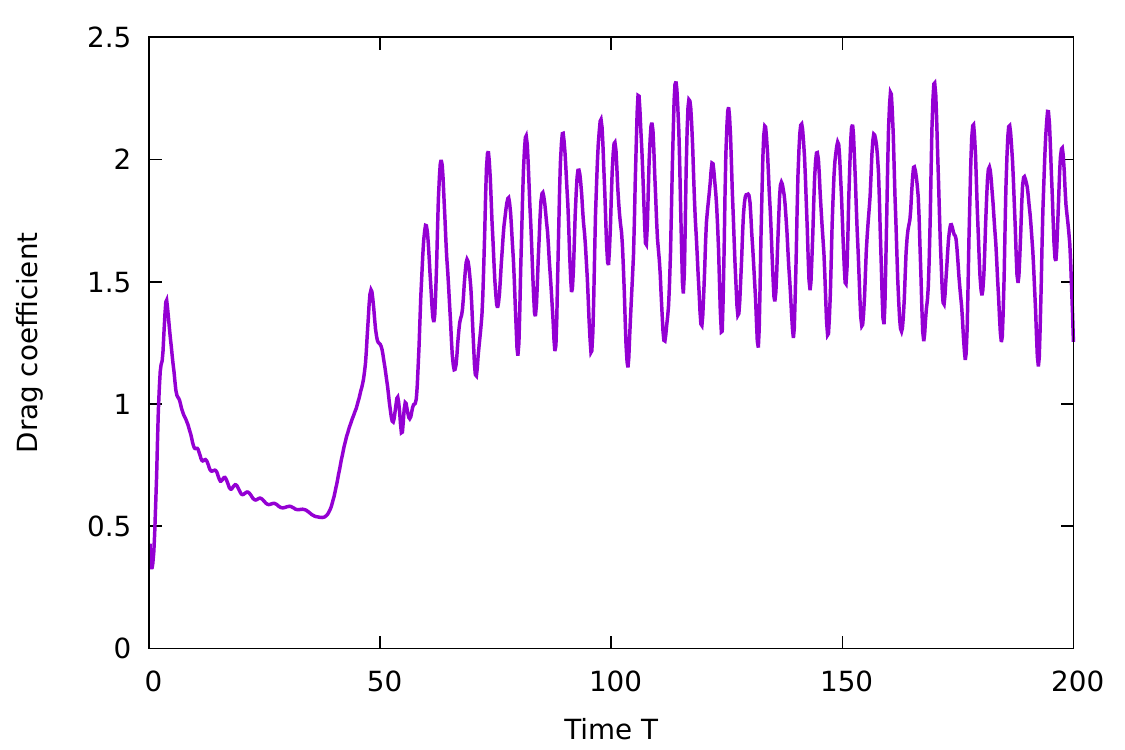}
\end{center}
\caption{Evolution of the drag coefficient of a circular cylinder for $\textrm{Re}=3000$ simulated on a grid with $4096\times 8192$ mesh points in the computational domain $\Omega=(-10,10)^2$.}
\label{fig:Drag_3000_Time_Long}
\end{figure*}

As previously mentioned, the flow remains symmetric at the beginning of the simulations for these flows around an impulsively started cylinder. By carrying the time integration over a much longer time
interval $T\in[0,200]$ instabilities due to round-off errors and to the nonlinearity of the system develop so that the flow becomes non symmetric for $T\ge 100$ at $\textrm{Re}=1000$ and $T\ge 50$ at
$\textrm{Re}=3000$ as it can be seen on Figures~\ref{fig:Drag_1000_Time_Long} and~\ref{fig:Drag_3000_Time_Long} representing the time history of the drag coefficient. After a transient period, an increase of the
drag coeffient is observed which stabilizes and oscillates around a mean value.

\subsection{Flow around moving bodies}
\label{Num_res_moving_bodies}
The purpose of this section is to show that the present numerical method is also able to simulate incompressible flows around moving bodies. 
Let us consider a cylinder which starts to move impulsively at $t = 0$ with the sinusoidal translational motion~
\begin{equation}
    \mathbf{u}_{\textrm{body}} (t) = \left( 2 \sin \left( \frac{t}{2} \right) ; 0 \right)
    \label{eq:TranslationalMotion}
\end{equation}
in a fluid initially at rest for Reynolds number $800$.
We suppose that the fluid is confined within a rectangular computational domain $\Omega = [-3;3] \times [-1;1]$ with no-slip boundary condition on $\partial \Omega^F$. 
The diameter of the cylinder is equal to $1$ and it is initially centered at the origin.
The boundary condition~\eqref{eq:TranslationalMotion} at the body surface $\partial \Omega^S$ is enforced through the non exhaustive following right hand side terms which vanish in the case of fixed obstacle : $u(\bs{\kappa}_{i,j}^{u,S})$, $u(\xi_{i,j-1/2}^S,y_{j-1/2})$, $\bs{u}_{\textrm{body}}(\bs{\kappa}_{i,j}^S)$ and so on. 
This terms are respectively taking into account in the convective terms, Laplacian operator and continuity equation. More details can be found in~\cite{BDJ1}. 

As the obstacle moves from one time step to another, we have to update matrices of the linear systems corresponding to Poisson and momentum equations at each iteration. 
Therefore, in such configuration, the direct solver requires much more CPU time compared to some iterative solver, which does not require a preprocessing step. 
For large problems, the faster solver we have found is an algebraic multigrid method (HYPRE BoomerAMG) implemented using the PETSc library~\cite{petsc-web-page,petsc-user-ref}.

\begin{figure*}
\begin{center}
  \includegraphics[width=0.8\textwidth]{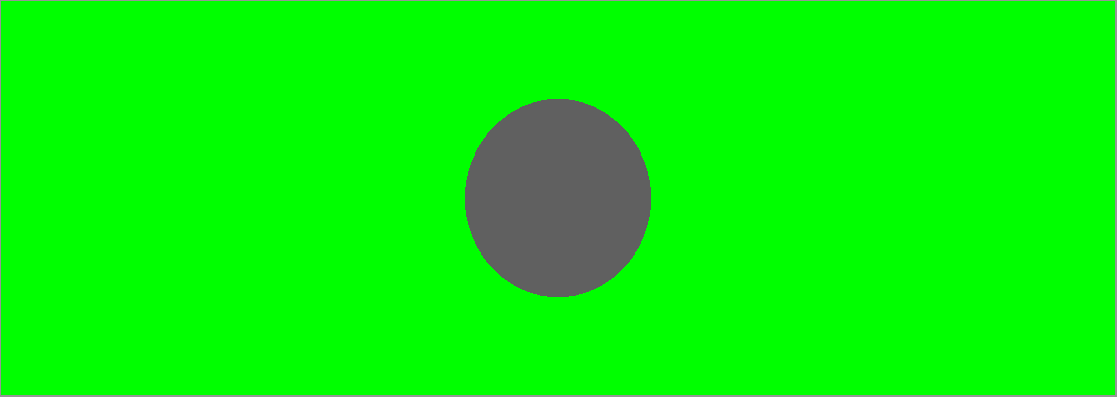} \\
        $t=0$ \\
  \includegraphics[width=0.8\textwidth]{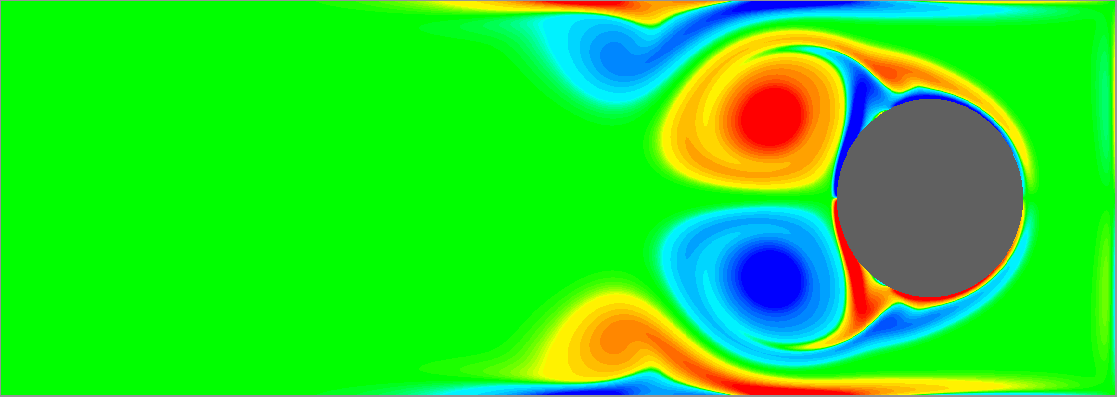} \\
        $t=\pi$ \\
  \includegraphics[width=0.8\textwidth]{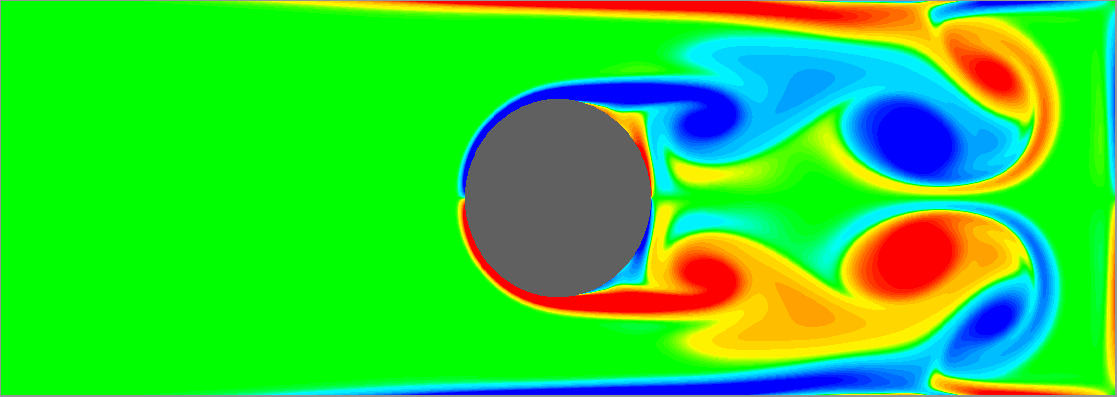} \\
          $t=2\pi$ \\
  \includegraphics[width=0.8\textwidth]{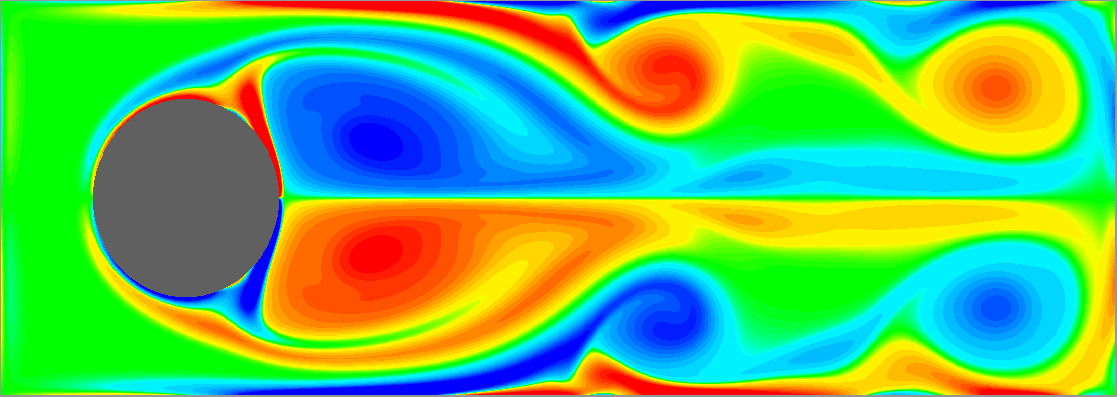} \\
          $t=3\pi$ \\
\end{center}
\caption{Flow around a moving circular cylinder at $\textrm{Re}=800$ in the computational domain $\Omega=(-3,3)\times(-1,1)$ discretized with $1200\times 400$ mesh points.}
\label{fig:Moving_Cylinder_800_1}
\end{figure*}

\begin{figure*}
\begin{center}
  \includegraphics[width=0.8\textwidth]{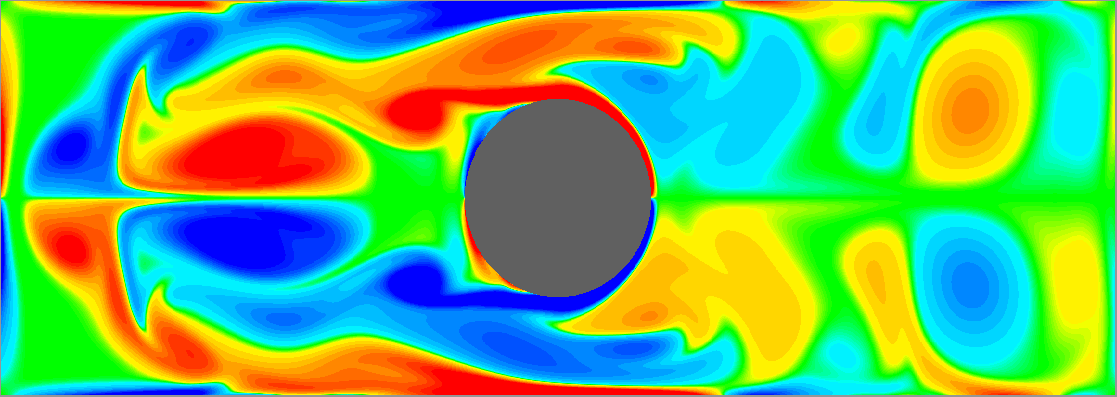} \\
            $T=4\pi$ \\
  \includegraphics[width=0.8\textwidth]{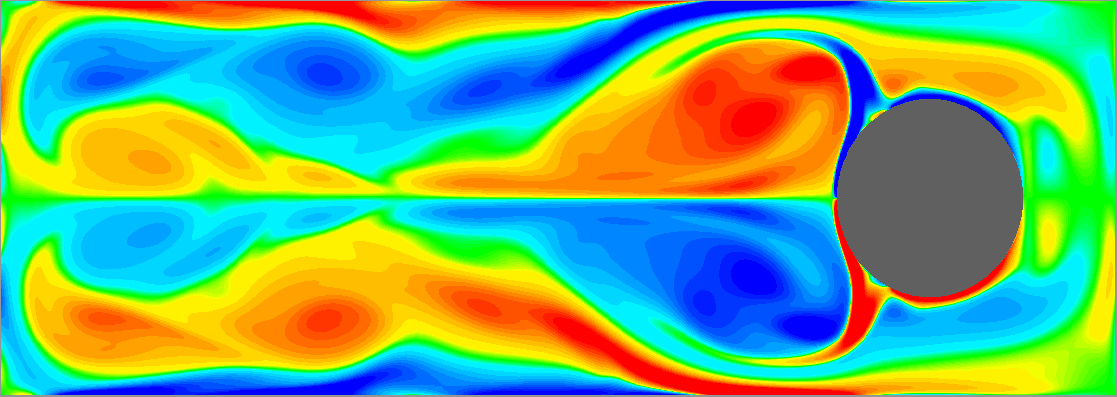} \\
            $T=5\pi$ \\
  \includegraphics[width=0.8\textwidth]{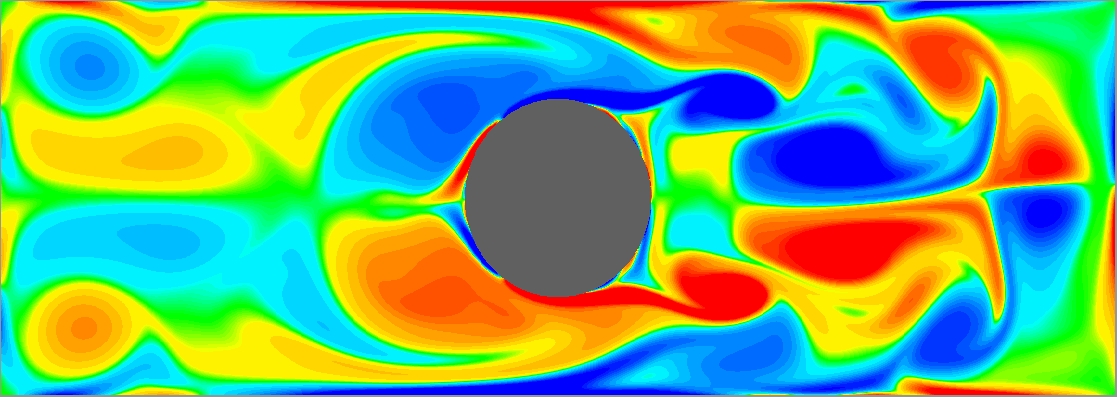} \\
            $T=6\pi$ \\
  \includegraphics[width=0.8\textwidth]{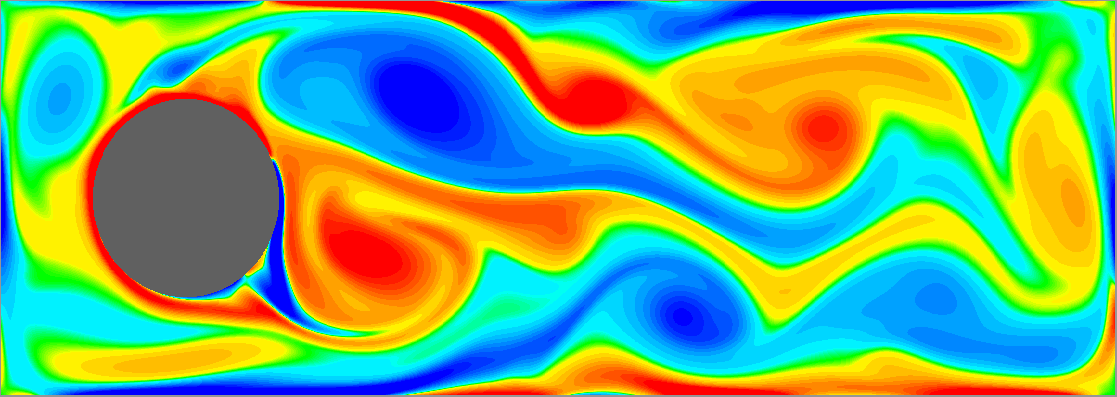} \\
            $T=7\pi$ \\
\end{center}
\caption{Flow around a moving circular cylinder at $\textrm{Re}=800$ in the computational domain $\Omega=(-3,3)\times(-1,1)$ discretized with $1200\times 400$ mesh points.}
\label{fig:Moving_Cylinder_800_2}
\end{figure*}

A constant mesh size $h=5\times 10^{-3}$ is used in both directions and the value of the time step, satisfying a CFL stability condition, is $10^{-3}$.
As shown in Figure~\ref{fig:Moving_Cylinder_800_1}, vortices interact with each other and also with the boundaries.
The flow remains perfectly symmetric until $t=6\pi$, thereafter the symmetry of the flow is lost due to rounding errors inherent in computer calculations (see Figure~\ref{fig:Moving_Cylinder_800_2}). 


\section{Conclusion}
\label{sec:conclusion}

We have presented a cut-cell method for the numerical solution of flows past obstacles. We have detailed the numerical method and the computational aspects for fixed obstacles, and shown numerical results for fixed and moving rigid bodies. 
The parallel version of the algorithm presented here allows computations of flows at Reynolds number up to 3000. The numerical tests confirm some results of the literature, a good agreement is observed with the numerical simulations in~\cite{KL1995}. We have also shown that the computation of the drag coefficient matches the theoretical square-root singularity predicted by~\cite{BY}.
The choice of the size of the box (compared with the grid size $h$) is one of the questions that we would like to investigate with this method in further works, we also would like to deal with rigid bodies following the fluid flow.

\section*{Acknowledgements}
This research was partially supported by the French Government Laboratory of Excellence initiative n$^\textrm{o}$ANR-10-LABX-0006,
the R\'egion Auvergne and the European Regional Development Fund. 
The numerical simulations have been performed on a DELL cluster with $32$ processors Xeon E2650v2 ($8$ cores), $1$ To of total
memory and an infiniband (FDR 56Gb/s) connecting network.
  
\bibliographystyle{plain}
\bibliography{bouchon_dubois_james.bib}

\end{document}